\documentclass{amsart}
\usepackage{amscd,amssymb}
\usepackage[mathcal]{euscript}

\numberwithin{equation}{section}

\newtheorem{main}{Theorem}

\newtheorem{theorem}[equation]{Theorem}
\newtheorem{lemma}[equation]{Lemma}
\newtheorem{proposition}[equation]{Proposition}

\newtheorem{corollary}[equation]{Corollary}
\theoremstyle{definition}
\newtheorem{definition}[equation]{Definition}
\newtheorem{remark}[equation]{Remark}

\DeclareMathOperator{\ho}{Ho}
\DeclareMathOperator{\colim}{colim}
\DeclareMathOperator{\invlim}{lim}
\DeclareMathOperator{\Ev}{Ev}

\DeclareMathOperator{\Arr}{Arr}
\DeclareMathOperator{\Sym}{Sym}

\DeclareMathOperator{\Hom}{Hom}
\DeclareMathOperator{\Map}{Map}
\DeclareMathOperator{\map}{map}

\DeclareMathOperator{\Stab}{Stab}

\newcommand{\cat}[1]{\mathcal{#1}}
\newcommand{\sset}{\mathbf{SSet}}
\renewcommand{\top}{\mathbf{Top}}

\newcommand{\symspec}{Sp^{\Sigma }}
\newcommand{\genBF}[2]{Sp^{\mathbb{N}}(#1,#2)}
\newcommand{\BF}{\genBF{\cat{C}}{G}}
\newcommand{\funcBF}[1]{Sp^{\mathbb{N}(#1)}}
\newcommand{\genspec}[2]{Sp^{\Sigma }(#1,#2)}
\newcommand{\spec}{\genspec{\cat{D}}{K}}
\newcommand{\modspec}{\genspec{\cat{C}}{K}}
\newcommand{\funcspec}{Sp^{\Sigma }(\Phi )}

\newcommand{\stensor}{\overline{\otimes }}
\newcommand{\boxprod}{\mathbin\square}
\renewcommand{\smash}{\wedge }
\renewcommand{\lim}{\invlim}
\newcommand{\cof}{\text{-cof}}
 
\newcommand{\inj}{\text{-inj}}

\newcommand{\rlp}{right lifting property with respect to }
\newcommand{\llp}{left lifting property with respect to }

\newcommand{\mathcolon}{\colon\,}
\newcommand{\ulp}{\textup{(}}
\newcommand{\urp}{\textup{)}}
\newcommand{\uc}{\textup{:}}

\newcommand{\ie}{\textit{i.e. }}

\begin{document}

\title{Spectra and symmetric spectra in general model categories}

\date{\today}
\author{Mark Hovey}
\address{Department of Mathematics \\ Wesleyan University 
\\ Middletown, CT 06459 \\ USA} 
\email{hovey@member.ams.org}



\begin{abstract}
We give two general constructions for the passage from unstable to
stable homotopy that apply to the known example of topological spaces,
but also to new situations, such as the $\mathbb{A}^{1}$-homotopy theory
of Morel-Voevodsky~\cite{morel-voevodsky, voevodsky}.  One is based on
the standard notion of spectra originated by
Boardman~\cite{boardman-vogt}.  Its input is a well-behaved model
category $\cat{C}$ and an endofunctor $G$, generalizing the suspension.
Its output is a model category $\BF $ on which $G$ is a Quillen
equivalence.  The second construction is based on symmetric
spectra~\cite{hovey-shipley-smith}, and applies to model categories
$\cat{C}$ with a compatible monoidal structure.  In this case, the
functor $G$ must be given by tensoring with a cofibrant object $K$.  The
output is again a model category $\spec $ where tensoring with $K$ is a
Quillen equivalence, but now $\spec $ is again a monoidal model
category.  We study general properties of these stabilizations; most
importantly, we give a sufficient condition for these two stabilizations
to be equivalent that applies both in the known case of topological
spaces and in the case of $\mathbf{A}^{1}$-homotopy theory.
\end{abstract}

\maketitle

\section*{Introduction}

The object of this paper is to give two very general constructions of
the passage from unstable homotopy theory to stable homotopy theory.
Since homotopy theory in some form appears in many different areas of
mathematics, this construction is useful beyond algebraic topology,
where these methods originated.  In particular, the two constructions we
give apply not only to the usual passage from unstable homotopy theory
of pointed topological spaces (or simplicial sets) to the stable
homotopy theory of spectra, but also to the passage from the unstable
$\mathbb{A}^{1}$-homotopy theory of Morel-Voevodsky~\cite{voevodsky,
morel-voevodsky} to the stable $\mathbb{A}^{1}$-homotopy theory.  This
example is obviously important, and the fact that it is an example of a
widely applicable theory of stabilization may come as a surprise to
readers of~\cite{jardine}, where specific properties of sheaves are
used.

Suppose, then, that we are given a (Quillen) model category $\cat{C}$
and a functor $G\mathcolon \cat{C}\xrightarrow{}\cat{C}$ that we would
like to invert, analogous to the suspension.  We will clearly need to
require that $\cat{C}$ be compatible with the model structure;
specifically, we require $G$ to be a left Quillen functor.  We will also
need some technical hypotheses on the model category $\cat{C}$, which
are complicated to state and to check, but which are satisfied in almost
all interesting examples, including $\mathbb{A}^{1}$-homotopy theory.
It is well-known what one should do to form the category $\BF $ of
spectra, as first written down for topological spaces
in~\cite{bousfield-friedlander}.  An object of $\BF $ is a sequence
$X_{n}$ of objects of $\cat{C}$ together with maps
$GX_{n}\xrightarrow{}X_{n+1}$, and a map $f\mathcolon X\xrightarrow{}Y$
is a sequence of maps $f_{n}\mathcolon X_{n}\xrightarrow{}Y_{n}$
compatible with the structure maps.  There is an obvious model
structure, called the \emph{projective model structure}, where the weak
equivalences are the maps $f\mathcolon X\xrightarrow{}Y$ such that
$f_{n}$ is a weak equivalence for all $n$.  It is not difficult to show
that this is a model structure and that there is a left Quillen functor
$G\mathcolon \BF \xrightarrow{}\BF $ extending $G$ on $\cat{C}$.  But,
just as in the topological case, $G$ will not be a Quillen equivalence.
So we must localize the projective model structure on $\BF $ to produce
the \emph{stable} model structure, with respect to which $G$ will be a
Quillen equivalence.  A new feature of this paper is that we are able to
construct the stable model structure with minimal hypotheses on
$\cat{C}$, using the localization results of
Hirschhorn~\cite{hirschhorn} (based on work of
Dror-Farjoun~\cite{dror}).  We must pay a price for this generality, of
course.  That price is that stable equivalences are not
stable homotopy isomorphisms, but instead are cohomology isomorphisms on all
cohomology theories, just as for symmetric
spectra~\cite{hovey-shipley-smith}.  If enough hypotheses are put on
$\cat{C}$ and $G$, then we show that stable equivalences are stable
homotopy isomorphisms.  Jardine~\cite{jardine} proves that stable
equivalences are stable homotopy isomorphisms in the stable
$\mathbb{A}^{1}$-homotopy theory, using the Nisnevitch descent theorem.  
His result does not follow from our general theorem for the
Morel-Voevodsky motivic model category, because the hypotheses we need
do not hold there, but Voevodsky (personal communication) has
constructed a simpler model category equivalent to the Morel-Voevodsky
one that does satify our hypotheses.  

As is well-known in algebraic topology, the category $\BF $ is not
sufficient to understand the smash product.  That is, if $\cat{C}$ is a
symmetric monoidal model category, and $G$ is the functor $X\mapsto
X\otimes K$ for some cofibrant object $K$ of $\cat{C}$, it almost never
happens that $\BF $ is symmetric monoidal.  We therefore need a
different construction in this case.  We define a category $\modspec $ just
as in symmetric spectra~\cite{hovey-shipley-smith}.  That is, an object
of $\modspec $ is a sequence $X_{n}$ of objects of $\cat{C}$ with an action
of the symmetric group $\Sigma _{n}$ on $X_{n}$.  In addition, we have
$\Sigma _{n}$-equivariant structure maps $X_{n}\otimes
K\xrightarrow{}X_{n+1}$, but we must further require that the iterated
structure maps $X_{n}\otimes K^{\otimes p}\xrightarrow{}X_{n+p}$ are
$\Sigma _{n}\times \Sigma _{p}$-equivariant, where $\Sigma _{p}$ acts on
$K^{\otimes p}$ by permuting the tensor factors.  It is once again
straightforward to construct the projective model structure on $\modspec $.
The same localization methods developed for $\BF $ apply again here to
give a stable model structure on which tensoring with $K$ is a Quillen
equivalence.  Once again, stable equivalences are cohomology
isomorphisms on all possible cohomology theories, but this time it is
very difficult to give a better description of stable equivalences even
in the case of simplicial symmetric spectra (but see~\cite{shipley-thh}
for the best such result I know).  We point out that our construction
gives a different construction of the stable model category of
simplicial symmetric spectra than the one appearing
in~\cite{hovey-shipley-smith}.

We now have competing stabilizations of $\cat{C}$ under the tensoring
with $K$ functor when $\cat{C}$ is symmetric monoidal.  Naturally, we
need to prove they are the same in an appropriate sense.  This was done
in the topological (actually, simplicial) case
in~\cite{hovey-shipley-smith} by constructing a functor $\BF
\xrightarrow{}\modspec $, where $K=S^{1}$ and $G$ is the tensor with
$S^{1}$ functor, and proving it is a Quillen equivalence.  We are unable
to generalize this argument.  Instead, following an idea of Hopkins, we
construct a zigzag of Quillen equivalences $\BF
\xrightarrow{}\cat{C}\xleftarrow{}\modspec $.  However, we need to require
that the cyclic permutation map on $K\otimes K\otimes K$ be homotopic to
the identity by an explicit homotopy for our construction to work.  This
hypothesis holds in the topological case with $K=S^{1}$ and in the
$\mathbb{A}^{1}$-local case with $K$ equal to either the simplicial
circle or the algebraic circle $\mathbb{A}^{1}-\{0 \}$.  This section of
the paper is by far the most delicate, and it is likely that we do not
have the best possible result.

We also investigate the properties of these two stabilization
constructions.  There are some obvious properties one would hope for of
a stabilization construction such as $\BF $ or $\modspec $.  First of all,
it should be functorial in the pair $(\cat{C},G)$.  We prove this for
both stabilization constructions; the most difficult point is defining
what one should mean by a map from $(\cat{C},G)$ to $(\cat{D},H)$.
Furthermore, it should be homotopy invariant.  That is, if the map
$(\cat{C},G)\xrightarrow{}(\cat{D},H)$ is a Quillen equivalence, the
induced map of stabilizations should also be a Quillen equivalence.  We
also prove this; one corollary is that the Quillen equivalence class of
$\modspec $ depends only on the homotopy type of $K$.  Finally, the
stabilization map $\cat{C}\xrightarrow{}\BF $ should be the initial map
to a model category $\cat{D}$ with an extension of $G$ to a Quillen
equivalence.  However, this last statement seems to be asking for too
much, because the category of model categories is itself something like
a model category.  This statement is analogous to asking for an initial
map in a model category from $X$ to a fibrant object, and such things do
not usually exist.  The best we can do is to say that if $G$ is already
a Quillen equivalence, then the map from $\cat{C}\xrightarrow{}\BF $ is
a Quillen equivalence.  This gives a weak form of uniqueness, and is the
basis for the comparison between $\BF $ and $\modspec $.  See also
see~\cite{schwede-shipley-unique} and~\cite{schwede-uniqueness} for
uniqueness results for the usual stable homotopy category.

We point out that this paper leaves some obvious questions open.  We do
not have a good characterization of stable equivalences or stable
fibrations in either spectra or symmetric spectra, in general, and we
are unable to prove that spectra or symmetric spectra are right proper.
We do have such characterizations for spectra when the original model
category $\cat{C}$ is sufficiently well-behaved, and the adjoint $U$ of
$G$ preserves sequential colimits.  These hypotheses include the cases
of ordinary simplicial spectra and spectra in a new motivic model
category of Voevodsky (but not the original Morel-Voevodsky motivic
model category).  We also prove that spectra are
right proper in this situation.  But we do not have a characterization
of stable equivalences of symmetric spectra even with these strong
assumptions.  Also, we have been unable to prove that symmetric spectra
satisfy the monoid axiom.  Without the monoid axiom, we do not get model
categories of monoids or of modules over an arbitrary monoid, though we
do get a model category of modules over a cofibrant monoid.  The
question of whether commutative monoids form a model category is even
more subtle and is not addressed in this paper.
See~\cite{mandell-may-shipley-schwede} for commutative monoids in
symmetric spectra of topological spaces.

There is a long history of work on stabilization, much of it not using
model categories.  As far as this author knows, Boardman was the first
to attempt to construct a good point-set version of spectra; his work
was never published (but see~\cite{boardman-vogt}), but it was the
standard for many years.  Generalizations of Boardman's construction
were given by Heller in several papers, including~\cite{heller-stable}
and~\cite{heller-completions}.   Heller has continued work on these
lines, most recently in~\cite{heller-polemics}.  The review of this
paper in Mathematical Reviews by Tony Elmendorf (MR98g:55021) captures
the response of many algebraic topologists to Heller's approach.  I
believe the central idea of Heller's approach is that the homotopy
theory associated to a model category $\cat{C}$ is the collection of all
possible homotopy categories of diagram categories $\text{ho }
\cat{C}^{I}$ and all functors between them.  With this definition, one
can then forget one had the model category in the first place, as Heller
does.  Unfortunately, the resulting complexity of definition is
overwhelming at present.

Of course, there has also been very successful work on stabilization by
May and coauthors, the two major milestones
being~\cite{lewis-may-steinberger}
and~\cite{elmendorf-kriz-mandell-may}.  At first glance, May's approach
seems wedded to the topological situation, relying as it does on
homeomorphisms $X_{n}\xrightarrow{}\Omega X_{n+1}$.  This is the reason
we have not tried to use it in this paper.  However, there has been
considerable recent work showing that this approach may be more flexible
than one might have expected.  I have
mentioned~\cite{mandell-may-shipley-schwede} above, but perhaps the most
ambitious attempt to generalize $S$-modules has been initiated by Mark
Johnson~\cite{johnson-mark-thesis}.  

Finally, we point out that Schwede~\cite{schwede} has shown that the
methods of Bousfield and Friedlander~\cite{bousfield-friedlander} apply
to certain more general model categories.  His model categories are
always simplicial and proper, and he is always inverting the ordinary
suspension functor.  Nevertheless, the paper~\cite{schwede} is the first
serious attempt to define a general stabilization functor of which the
author is aware.

This paper is organized as follows.  We begin by defining the category
$\BF $ and the associated strict model structure in
Section~\ref{sec-G-spectra}.  Then there is the brief
Section~\ref{sec-bousfield} recalling Hirschhorn's approach to
localization of model categories.  We then construct the stable model
structure modulo certain technical lemmas in Section~\ref{sec-stable}.
The technical lemmas we need assert that if a model category $\cat{C}$
is left proper cellular, then so is the strict model structure on $\BF
$, and therefore we can apply the localization technology of Hirschhorn.
We prove these technical lemmas, and the analogous lemmas for the strict
model structure on $\spec $, in an Appendix.  In
Section~\ref{sec-fin-gen}, we study the simplifications that arise when
the adjoint $U$ of $G$ preserves sequential colimits and $\cat{C}$ is
sufficiently well-behaved.  We characterize stable equivalences as the
appropriate generalization of stable homotopy isomorphisms in this case,
and we show the stable model structure is right proper, giving a
description of the stable fibrations as well.  In
Section~\ref{sec-funct}, we prove the functoriality, homotopy
invariance, and homotopy idempotence of the construction
$(\cat{C},G)\mapsto \BF $.  We also investigate monoidal structure.

Section~\ref{sec-symm-spec} begins the second part of the paper, about
symmetric spectra.  Since we have developed all the necessary techniques
in the first part, the proofs in this part are more concise.  In
Section~\ref{sec-symm-spec} we discuss the category of symmetric
spectra.  In Section~\ref{sec-symmetric-model} we construct the
projective and stable model structures on symmetric spectra, and in
Section~\ref{sec-symmetric-functor}, we discuss some properties of
symmetric spectra.  This includes functoriality, homotopy invariance,
and homotopy idempotence of the stable model structure.  We conclude the
paper in Section~\ref{sec-comparison} by constructing the chain of
Quillen equivalences betwee $\BF $ and $\modspec $, under the cyclic
permutation hypothesis mentioned above.  Finally, as stated previously,
there is an Appendix verifying that the techniques of Hirschhorn can be
applied to the projective model structures on $\BF $ and $\modspec $.

Obviously, considerable familiarity with model categories will be
necessary to understand this paper.  The original reference
is~\cite{quillen-htpy}, but a better introductory reference
is~\cite{dwyer-spalinski}.  More in depth references
include~\cite{hirschhorn}, ~\cite{hovey-model}, and~\cite{kan-model}.
In particular, we rely heavily on the localization technology
in~\cite{hirschhorn}.  

The author would like to thank Dan Dugger, Phil Hirschhorn, Mike
Hopkins, Dan Kan, Stefan Schwede, Brooke Shipley, Jeff Smith, Markus
Spitzweck, and Vladimir Voevodsky for helpful conversations about this
paper.  In particular, to the author's knowledge, it is Jeff Smith's
vision that one should be able to stabilize an arbitrary model category,
a vision that could not be carried out without Phil Hirschhorn's
devotion to getting the localization theory of model categories right.
The idea of using almost finitely generated model categories in
Section~\ref{sec-fin-gen} is due to Voevodsky, and the idea of using
bispectra to compare symmetric spectra with ordinary spectra (see
Section~\ref{sec-comparison}) is due to Hopkins.

\section{Spectra}\label{sec-G-spectra}

In this section and throughout the paper, $\cat{C}$ will be a model
category and $G\mathcolon \cat{C}\xrightarrow{}\cat{C}$ will be a
left Quillen endofunctor of $\cat{C}$ with right adjoint $U$.  In this
section, we define the category $\BF $ of spectra and construct its
strict model structure.  

The following definition is a straightforward generalization of the usual
notion of spectra~\cite{bousfield-friedlander}.  

\begin{definition}\label{defn-G-spectra}
Suppose $G$ is a left Quillen endofunctor of a model category $\cat{C}$.
Define $\BF $, the category of \emph{spectra}, as follows.  A
\emph{spectrum} $X$ is a sequence $X_{0},X_{1},\dots ,X_{n},\dots $ of
objects of $\cat{C}$ together with structure maps $\sigma \mathcolon
GX_{n}\xrightarrow{}X_{n+1}$ for all $n$.  A \emph{map of spectra} from
$X$ to $Y$ is a collection of maps $f_{n}\mathcolon
X_{n}\xrightarrow{}Y_{n}$ commuting with the structure maps, so that the
diagram
\[
\begin{CD}
GX_{n} @>\sigma _{X}>> X_{n+1} \\
@VGf_{n}VV @VVf_{n+1}V \\
GY_{n} @>>\sigma _{Y}> Y_{n+1}
\end{CD}
\]
is commutative for all $n$.  
\end{definition}

Note that if $\cat{C}$ is either the model category of pointed
simplicial sets or the model category of pointed topological spaces, and
$G$ is the suspension functor given by smashing with the circle $S^{1}$,
then $\BF $ is the Bousfield-Friedlander category of
spectra~\cite{bousfield-friedlander}.

\begin{definition}\label{defn-ev}
Given $n\geq 0$, the \emph{evaluation functor} $\Ev _{n}\mathcolon
\BF \xrightarrow{}\cat{C}$ takes $X$ to $X_{n}$.  The evaluation
functor has a left adjoint $F_{n}\mathcolon
\cat{C}\xrightarrow{}\BF $ defined by $(F_{n}X)_{m}=G^{m-n}X$ if
$m\geq n$ and $(F_{n}X)_{m}=0$ otherwise, where $0$ is the initial object
of $\cat{C}$.  The structure maps are the obvious ones.  
\end{definition}

Note that $F_{0}$ is an full and faithful embedding of the category
$\cat{C}$ into $\BF $.

Limits and colimits in the category of spectra are taken objectwise.  

\begin{lemma}\label{lem-spectra-complete}
The category of spectra is bicomplete.  
\end{lemma}

\begin{proof}
Given a functor $X\mathcolon \cat{I}\xrightarrow{}\BF $, we
define $(\colim X)_{n}=\colim \Ev_{n}\circ X$ and $(\lim X)_{n}=\lim \Ev
_{n}\circ X$.  Since $G$ is a left adjoint, it preserves colimits.  The
structure maps of the colimit are then the composites 
\[
G(\colim \Ev _{n}\circ X)\cong \colim (G\circ \Ev _{n}\circ X)
\xrightarrow{\colim \sigma \circ X}\colim \Ev _{n+1}\circ X  .
\]
Although $G$ does not preserve limits, there is still a natural map
$G(\lim Y) \xrightarrow{} \lim GY$ for any functor $Y \mathcolon \cat{I}
\xrightarrow{} \cat{C}$.  Then the structure maps of the limit are the
composites 
\[
G(\lim \Ev _{n}\circ X)\xrightarrow{}\lim (G\circ \Ev _{n}\circ X)
\xrightarrow{\lim \sigma \circ X} \lim \Ev _{n+1} \circ X . \qed
\]
\renewcommand{\qed}{}
\end{proof}

\begin{remark}\label{rem-Rn}
Note that the evaluation functor $\Ev _{n}\mathcolon
\BF \xrightarrow{}\cat{C}$ preserves colimits, so should have a
right adjoint $R_{n}\mathcolon \cat{C}\xrightarrow{}\BF $.  We
define $(R_{n}X)_{m}=U^{n-m}X$ if $m\leq n$, and $(R_{n}X)_{m}=1$ if
$m>n$.  The structure map $GU^{n-m}X\xrightarrow{}U^{n-m-1}X$ is adjoint
to the identity map of $U^{n-m}X$ when $m<n$.  We leave it to the reader
to verify that $R_{n}$ is the right adjoint of $\Ev _{n}$.
\end{remark}

We now show that the functors $G$ and $U$ extend to functors on
$\BF $.  

\begin{lemma}\label{lem-extending-functors}
Suppose $F\mathcolon \cat{C}\xrightarrow{}\cat{C}$ is a functor and
$\tau \mathcolon GF\xrightarrow{}FG$ is a natural transformation.  Then
there is an induced functor $F\mathcolon \BF 
\xrightarrow{}\BF $, called the \emph{prolongation} of $F$. 
\end{lemma}

\begin{proof}
Given $X\in \BF $, we define $(FX)_{n}=FX_{n}$, with structure
maps $GFX_{n}\xrightarrow{\tau }FGX_{n}\xrightarrow{F\sigma }FX_{n+1}$.
Given a map $f\mathcolon X\xrightarrow{}Y$, we define $Ff$ by
$(Ff)_{n}=Ff_{n}$.  
\end{proof}

Note that the prolongation of $F$ depends on the choice of natural
transformation $\tau $.  Usually there is an obvious choice of $\tau $. 

\begin{corollary}\label{cor-G-extends}
The functors $G$ and $U$ prolong to adjoint functors $G$ and $U$ on
$\BF $.
\end{corollary}

\begin{proof}
To prolong $G$, take $\tau $ to be the identity in
Lemma~\ref{lem-extending-functors}.  To prolong its right adjoint $U$,
take $\tau $ to be the composite $GUX\xrightarrow{}X\xrightarrow{}UGX$
of the counit and unit of the adjunction.  It is a somewhat involved
exercise in adjoint functors to verify that the resulting prolongations
$G'$ and $U'$ are still adjoint to each other.  Like all the exercises in
adjoint functors in this paper, the simplest way to proceed is to write
the adjoint $Y\xrightarrow{}UZ$ of a map $f\mathcolon GY\xrightarrow{}Z$
as the composition
\[
Y \xrightarrow{\eta } UGY \xrightarrow{Uf} UZ
\]
and to use the standard properties of the unit and counit of an
adjunction.  
\end{proof}

The following remark is critically important to the understanding of our
approach to spectra.  

\begin{remark}\label{rem-subtle}
The definition we have just given of the prolongation of $G$ to an
endofunctor of $\BF $ is the only possible definition under our very
general hypotheses.  However, this definition \textbf{does not
generalize the definition of the suspension} when $\cat{C}$ is the
category of pointed topological spaces and $GX=X\wedge S^{1}$.  Indeed,
recall from~\cite{bousfield-friedlander} that the suspension of a
spectrum $X$ in this case is defined by $(X\wedge S^{1})_{n}=X_{n}\wedge
S^{1}$, with structure map given by 
\[
X_{n}\wedge S^{1}\wedge S^{1} \xrightarrow{1\wedge T} X_{n} \wedge S^{1} \wedge
S^{1} \xrightarrow{\sigma \wedge 1} X_{n+1}\wedge S^{1},
\]
where $T$ is the twist isomorphism.  On the other hand, if we apply our
definition of the prolongation of $G$ above, we get a functor $X\mapsto
X\otimes S^{1}$ defined by $(X\otimes S^{1})_{n}=X_{n}\otimes S^{1}$
with structure map
\[
X_{n}\wedge S^{1} \wedge S^{1} \xrightarrow{\sigma \wedge 1}
X_{n+1}\wedge S^{1}.  
\]
Said another way, Bousfield and Friedlander choose the natural
transformation $\tau \mathcolon X\wedge S^{1}\wedge
S^{1}\xrightarrow{}X\wedge S^{1}\wedge S^{1}$ to be $1\wedge T$, while
we are taking it to be the identity.  This is a crucial and subtle
difference whose ramifications we will study in
Section~\ref{sec-comparison}.  
\end{remark}

We now show that $\BF $ inherits a model structure from $\cat{C}$,
called the \emph{projective model structure}.  The functor $G\mathcolon
\BF \xrightarrow{}\BF $ will be a left Quillen functor with respect to
the projective model structure, but it will not be a Quillen
equivalence.  Our approach to the projective model structure owes much
to~\cite{bousfield-friedlander}
and~\cite[Section~5.1]{hovey-shipley-smith}.  At this point, we will
slip into the standard model category terminology and notation, all of
which can be found in~\cite{hovey-model}, mostly in Section~2.1.  

\begin{definition}\label{defn-level}
A map $f\in \BF $ is a \emph{level equivalence} if each map
$f_{n}$ is a weak equivalence in $\cat{C}$.  Similarly, $f$ is a
\emph{level fibration} (resp.\emph{level cofibration}, \emph{level
trivial fibration}, \emph{level trivial cofibration}) if each map
$f_{n}$ is a fibration (resp. cofibration, trivial fibration, trivial
cofibration) in $\cat{C}$.  The map $f$ is a \emph{projective
cofibration} if $f$ has the \llp every level trivial fibration.  
\end{definition}

Note that level equivalences satisfy the two out of three property, and
each of the classes defined above is closed under retracts.  Thus we
should be able to construct a model structure using these classes.  To
do so, we need the small object argument, and hence we assume that
$\cat{C}$ is cofibrantly generated (see~\cite[Section~2.1]{hovey-model}
for a discussion of cofibrantly generated model categories).

\begin{definition}\label{defn-IG}
Suppose $\cat{C}$ is a cofibrantly generated model category with
generating cofibrations $I$ and generating trivial cofibrations $J$.
Suppose $G$ is a left Quillen endofunctor of $\cat{C}$, and form the
category of spectra $\BF $.  Define sets of maps in $\BF $ by
$I_{G}=\bigcup _{n} F_{n}I$ and $J_{G}=\bigcup _{n}F_{n}J$.
\end{definition}

The sets $I_{G}$ and $J_{G}$ will be the generating cofibrations and
trivial cofibrations for a model structure on $\BF $.  There is a
standard method for proving this, based on the small object
argument~\cite[Theorem~2.1.14]{hovey-model}.  The first step is to show
that the domains of $I_{G}$ and $J_{G}$ are small, in the sense
of~\cite[Definition~2.1.3]{hovey-model}. 

\begin{proposition}\label{prop-smallness-in-G-spec}
Suppose $A$ is small relative to the cofibrations \ulp resp. trivial
cofibrations\urp{} in $\cat{C}$, and $n\geq 0$.  Then $F_{n}A$ is small
relative to the level cofibrations \ulp resp. level trivial
cofibrations\urp{} in $\BF $.
\end{proposition}

\begin{proof}
The main point is that $\Ev _{n}$ commutes with colimits.  We leave the
remainder of the proof to the reader.  
\end{proof}

To apply this to the domains of $I_{G}$, we need to know that the maps
of $I_{G}\cof $ are level cofibrations.  Recall the right adjoint
$R_{n}$ of $\Ev _{n}$ constructed in Remark~\ref{rem-Rn}.  

\begin{lemma}\label{lem-char-level-cofib}
A map in $\BF $ is a level cofibration \ulp resp. level trivial
cofibration\urp{} if and only if it has the \llp $R_{n}g$ for all $n\geq
0$ and all trivial fibrations \ulp resp. fibrations\urp{} $g$ in
$\cat{C}$.
\end{lemma}

\begin{proof}
By adjunction, a map $f$ has the \llp $R_{n}g$ if and only if $\Ev
_{n}f$ has the \llp $g$.  Since a map is a cofibration (resp. trivial
cofibration) in $\cat{C}$ if and only if it has the \llp all trivial
fibrations (resp. fibrations), the lemma follows.
\end{proof}

\begin{proposition}\label{prop-proj-cofs-are-level}
Every map in $I_{G}\cof $ is a level cofibration.  Every map in
$J_{G}\cof $ is a level trivial cofibration.  
\end{proposition}

\begin{proof}
Since $G$ is a left Quillen functor, every map in $I_{G}$ is a level
cofibration.  By Lemma~\ref{lem-char-level-cofib}, this means that
$R_{n}g\in I_{G}\inj $ for all $n\geq 0$ and all trivial fibrations $g$.
Since a map in $I_{G}\cof $ has the \llp every map in $I_{G}\inj $, in
particular it has the \llp $R_{n}g$.  Another application of
Lemma~\ref{lem-char-level-cofib} completes the proof for $I_{G}\cof $.
The proof for $J_{G}\cof $ is similar.
\end{proof}

Proposition~\ref{prop-smallness-in-G-spec} and
Proposition~\ref{prop-proj-cofs-are-level} immediately imply the
following corollary.  

\begin{corollary}\label{cor-IG-small}
The domains of $I_{G}$ are small relative to $I_{G}\cof $.  The domains
of $J_{G}$ are small relative to $J_{G}\cof $.
\end{corollary}

\begin{theorem}\label{thm-strict}
Suppose $\cat{C}$ is cofibrantly generated.  Then the projective
cofibrations, the level fibrations, and the level equivalences define a
cofibrantly generated model structure on $\BF $, with generating
cofibrations $I_{G}$ and generating trivial cofibrations $J_{G}$.  We
call this the \emph{projective model structure}.  The projective model
structure is left proper \ulp resp. right proper, proper\urp{} if
$\cat{C}$ is left proper \ulp resp. right proper, proper\urp .
\end{theorem}

Note that if $\cat{C}$ is either the model category of pointed
simplicial sets or pointed topological spaces, and $G$ is the suspension
functor, the projective model structure on $\cat{C}_{G}$ is the strict
model structure on the Bousfield-Friedlander category of
spectra~\cite{bousfield-friedlander}.


\begin{proof}
The retract and two out of three axioms are immediate, as is the lifting
axiom for a projective cofibration and a level trivial fibration.  By
adjointness, a map is a level trivial fibration if and only if it is in
$I_{G}\inj $.  Hence a map is a projective cofibration if and only if it
is in $I_{G}\cof $.  The small object
argument~\cite[Theorem~2.1.14]{hovey-model} applied to $I_{G}$ then
produces a functorial factorization into a projective cofibration
followed by a level trivial fibration.

Adjointness implies that a map is a level fibration if and only if it is
in $J_{G}\inj $.  We have already seen in
Proposition~\ref{prop-proj-cofs-are-level} that the maps in $J_{G}\cof $
are level equivalences, and they are projective cofibrations since they
have the \llp all level fibrations, and in particular level trivial
fibrations.  Hence the small object argument applied to $J_{G}$ produces
a functorial factorization into a projective cofibration and level
equivalence followed by a level fibration.  

Conversely, we claim that any projective cofibration and level
equivalence $f$ is in $J_{G}\cof $, and hence has the \llp level
fibrations.  To see this, write $f=pi$ where $i$ is in $J_{G}\cof $ and
$p$ is in $J_{G}\inj $.  Then $p$ is a level fibration.  Since $f$ and
$i$ are both level equivalences, so is $p$.  Thus $f$ has the \llp $p$,
and so $f$ is a retract of $i$ by the retract
argument~\cite[Lemma~1.1.9]{hovey-model}.  In particular $f\in J_{G}\cof
$.

Since colimits and limits in $\BF $ are taken levelwise, and since every
projective cofibration is in particular a level cofibration, the
statements about properness are immediate.
\end{proof}

We also characterize the projective cofibrations.  We denote the pushout
of two maps $A\xrightarrow{}B$ and $A\xrightarrow{}C$ by $B\amalg
_{A}C$.

\begin{proposition}\label{prop-char-of-proj-cofib}
A map $i\mathcolon A\xrightarrow{}B$ is a projective \ulp trivial\urp{}
cofibration if and only if the induced maps $A_{0}\xrightarrow{}B_{0}$
and $A_{n}\amalg _{GA_{n-1}}GB_{n-1}\xrightarrow{}B_{n}$ for $n\geq 1$
are \ulp trivial\urp{} cofibrations.
\end{proposition}

\begin{proof}
We only prove the cofibration case, leaving the similar trivial
cofibration case to the reader.  First suppose $i\mathcolon
A\xrightarrow{}B$ is a projective cofibration.  We have already seen in
Proposition~\ref{prop-proj-cofs-are-level} that
$A_{0}\xrightarrow{}B_{0}$ is a cofibration.  Suppose $p\mathcolon
X\xrightarrow{}Y$ is a trivial fibration in $\cat{C}$, and suppose we
have a commutative diagram 
\[
\begin{CD}
A_{n}\amalg _{GA_{n-1}}GB_{n-1} @>>> X \\
@VVV @VVpV \\
B_{n} @>>> Y
\end{CD}
\]
We must construct a lift in this diagram.  By adjointness, it suffices
to construct a lift in the induced diagram 
\[
\begin{CD}
A @>>> R_{n}X \\
@ViVV @VVV \\
B @>>> R_{n}Y\times _{R_{n-1}UY} R_{n-1}UX
\end{CD}
\]
where $R_{n}$ is the right adjoint of $\Ev _{n}$.  Using the description
of $R_{n}$ given in Remark~\ref{rem-Rn}, one can easily check that the
map $R_{n}X\xrightarrow{}R_{n}Y\times _{R_{n-1}UY} R_{n-1}UX$ is a level
trivial fibration, so a lift exists.  

Conversely, suppose that the map $i$ satisfies the conditions in the
statement of the proposition.  Suppose $p\mathcolon X\xrightarrow{}Y$ is
a level trivial fibration in $\BF $, and suppose the diagram
\[
\begin{CD}
A @>f>> X \\
@ViVV @VVpV \\
B @>>g> Y
\end{CD}
\]
commutes.  We construct a lift $h_{n}\mathcolon
B_{n}\xrightarrow{}X_{n}$, compatible with the structure maps, by
induction on $n$.  There is no difficulty defining $h_{0}$, since
$i_{0}$ has the \llp the trivial fibration $p_{0}$.  Suppose we have
defined $h_{j}$ for $j<n$.  Then by lifting in the induced diagram
\[
\begin{CD}
A_{n}\amalg _{GA_{n-1}}GB_{n-1} @>(f_{n},\sigma \circ Gh_{n-1})>> X_{n} \\
@VVV @VVp_{n}V \\
B_{n} @>>g_{n}> Y_{n}
\end{CD}
\]
we find the required map $h_{n}\mathcolon B_{n}\xrightarrow{}X_{n}$.  
\end{proof}

Finally, we point out that the prolongation of $G$ is still a Quillen
functor.

\begin{proposition}\label{prop-G-Quillen-functor}
Give $\BF $ the projective model structure.  Then the
prolongation $G\mathcolon \BF \xrightarrow{}\BF $ of $G$
is a Quillen functor.  Furthermore, the functor $F_{n}\mathcolon
\cat{C}\xrightarrow{}\BF $ is a Quillen functor.
\end{proposition}

\begin{proof}
The functor $\Ev _{n}$ obviously takes level fibrations to fibrations
and level trivial fibrations to trivial fibrations.  Hence $\Ev _{n}$ is
a right Quillen functor, and so its left adjoint $F_{n}$ is a left
Quillen functor.  Similarly, the prolongation of $U$ to a functor
$U\mathcolon \BF \xrightarrow{}\BF $ preserves level fibrations and
level trivial fibrations, so its left adjoint $G$ is a Quillen functor.
\end{proof}

\section{Bousfield localization}\label{sec-bousfield}

We will define the stable model structure on $\BF$ in
Section~\ref{sec-stable} as a Bousfield localization of the projective
model structure on $\BF $.  In this section we recall the theory
of Bousfield localization of model categories from~\cite{hirschhorn}.

To do so, we need some preliminary remarks related to function
complexes.  Details can be found in~\cite[Chapter 5]{hovey-model},
\cite{kan-model}, and~\cite[Chapter 19]{hirschhorn}.  Given an object
$A$ in a model category $\cat{C}$, there is a functorial cosimplicial
resolution of $A$ induced by the functorial factorizations of $\cat{C}$.
By mapping out of this cosimplicial resolution we get a simplicial set
$\Map _{\ell }(A,X)$.  Similarly, there is a functorial simplicial
resolution of $X$, and by mapping into it we get a simplicial set $\Map
_{r}(A,X)$.  One should think of these as replacements for the
simplicial structure present in a simplicial model category.  Let us
define the homotopy function complex $\map (A,X)=\Map _{r}(QA,RX)$.
Then $\map (A,X)$ is canonically isomorphic in the homotopy category of
simplicial sets to $\Map_{\ell }(QA,RX)$, and makes $\ho \cat{C}$ an
enriched category over the homotopy category $\ho \sset $ of simplicial
sets.  In fact, $\ho \cat{C}$ is naturally tensored and cotensored over
$\ho \sset $, as well as enriched over it.  In particular, if $G$ is an
arbitrary left Quillen functor between model categories with right
adjoint $U$, we have $\map ((LG)X,Y)\cong \map (X,(RU)Y)$ in $\ho \sset
$, where $(LG)X=GQX$ is the total left derived functor of $G$ and
$(RU)Y=URY$ is the total right derived functor of $U$.

\begin{definition}\label{defn-S-local}
Suppose we have a set $\cat{S}$ of maps in a model category $\cat{C}$.
\begin{enumerate}
\item A \emph{$\cat{S}$-local object} of $\cat{C}$ is a fibrant object
$W$ such that, for every $f\mathcolon A\xrightarrow{}B$ in $\cat{S}$,
the induced map $\map(B,W)\xrightarrow{}\map(A,W)$ is a weak equivalence
of simplicial sets.
\item A \emph{$\cat{S}$-local equivalence} is a map $g\mathcolon
A\xrightarrow{}B$ in $\cat{C}$ such that the induced map $\map
(B,W)\xrightarrow{}\map (A,W)$ is a weak equivalence of simplicial sets
for all $\cat{S}$-local objects $W$.
\end{enumerate}
\end{definition}

Note that $\cat{S}$-local equivalences between $\cat{S}$-local objects
are in fact weak equivalences.

Then the main theorem of~\cite{hirschhorn} is the following.  We will
define cellular model categories, a special class of cofibrantly
generated model categories, in the Appendix.

\begin{theorem}\label{thm-localization}
Suppose $\cat{S}$ is a set of maps in a left proper cellular model
category $\cat{C}$.  Then there is a left proper cellular model
structure on $\cat{C}$ where the weak equivalences are the
$\cat{S}$-local equivalences and the cofibrations remain unchanged.  The
$\cat{S}$-local objects are the fibrant objects in this model structure.
We denote this new model category by $L_{\cat{S}}\cat{C}$ and refer to
it as the \emph{Bousfield localization} of $\cat{C}$ with respect to
$\cat{S}$.  Left Quillen functors from $L_{\cat{S}}\cat{C}$ to $\cat{D}$
are in one to one correspondence with left Quillen functors $F\mathcolon
\cat{C}\xrightarrow{}\cat{D}$ such that $F(Qf)$ is a weak equivalence
for all $f\in \cat{S}$.
\end{theorem}


We will also need the following fact about localizations, which is
implicit in~\cite[Chapter 4]{hirschhorn}.  

\begin{proposition}\label{prop-local-quillen}
Suppose $\cat{C}$ and $\cat{D}$ are left proper cellular model
categories, $\cat{S}$ is a set of maps in $\cat{C}$, and $\cat{T}$ is a
set of maps in $\cat{D}$.  Suppose $F\mathcolon
\cat{C}\xrightarrow{}\cat{D}$ is a Quillen equivalence with right
adjoint $U$, and suppose $F(Qf)$ is a $\cat{T}$-local equivalence for
all $f\in \cat{S}$.  Then $F$ induces a Quillen \emph{equivalence}
$F\mathcolon L_{\cat{S}}\cat{C}\xrightarrow{}L_{\cat{T}}\cat{D}$ if and only
if, for every $\cat{S}$-local $X\in \cat{C}$, there is a $\cat{T}$-local
$Y$ in $\cat{D}$ such that $X$ is weakly equivalent in $\cat{C}$ to
$UY$.  This condition will hold if, for all fibrant $Y$ in $\cat{D}$
such that $UY$ is $\cat{S}$-local, $Y$ is $\cat{T}$-local.
\end{proposition}

\begin{proof}
Suppose first that $F$ does induce a Quillen equivalence on the
localizations, and suppose that $X$ is $\cat{S}$-local.  Then $QX$ is
also $\cat{S}$-local.  Let $L_{\cat{T}}$ denote a fibrant replacement
functor in $L_{\cat{T}}\cat{D}$.  Then, because $F$ is a Quillen
equivalence on the localizations, the map
$QX\xrightarrow{}UL_{\cat{T}}FQX$ is a weak equivalence in
$L_{\cat{S}}\cat{C}$ (see~\cite[Section~1.3.3]{hovey-model}).  But both
$QX$ and $UL_{\cat{T}}FQX$ are $\cat{S}$-local, so
$QX\xrightarrow{}UL_{\cat{T}}FQX$ is a weak equivalence in $\cat{C}$.
Hence $X$ is weakly equivalent in $\cat{C}$ to $UY$, where $Y$ is the
$\cat{T}$-local object $L_{\cat{T}}FQX$.

Before proving the converse, note that, since $F$ is a Quillen
equivalence before localizing, the map $FQUX\xrightarrow{}X$ is a weak
equivalence for all fibrant $X$.  Since the functor $Q$ does not change
upon localization, and every $\cat{T}$-local object of $\cat{D}$ is in
particular fibrant, this condition still holds after localization.  Thus
$F$ is a Quillen equivalence after localization if and only if $F$
reflects local equivalences between cofibrant objects,
by~\cite[Corollary~1.3.16]{hovey-model}.  

Suppose $f\mathcolon A\xrightarrow{}B$ is a map between cofibrant
objects such that $Ff$ is a $\cat{T}$-local equivalence.  We must show
that $\map (f,X)$ is a weak equivalence for all $\cat{S}$-local $X$.
Adjointness implies that $\map (f,UY)$ is a weak equivalence for all
$\cat{T}$-local $Y$, and our condition then guarantees that this is
enough to conclude that $\map (f,X)$ is a weak equivalence for all
$\cat{S}$-local $X$.

We still need to prove the last statement of the proposition.  So
suppose $X$ is $\cat{S}$-local.  Then $QX$ is also $\cat{S}$-local, and,
in $\cat{C}$, we have a weak equivalence $QX\xrightarrow{}URFQX$.  Our
assumption then guarantees that $Y=RFQX$ is $\cat{T}$-local, and $X$ is
in indeed weakly equivalent to $UY$.
\end{proof}

The fibrations in $L_{\cat{S}}\cat{C}$ are not completely
understood~\cite[Section~3.6]{hirschhorn}.  The $\cat{S}$-local
fibrations between $\cat{S}$-local fibrant objects are just the usual
fibrations.  In case both $\cat{C}$ and $L_{\cat{S}}\cat{C}$ are right
proper, there is a characterization of the $\cat{S}$-local fibrations in
terms of homotopy pullbacks analogous to the characterization of stable
fibrations of spectra in~\cite{bousfield-friedlander}.  However,
$L_{\cat{S}}\cat{C}$ need not be right proper even if $\cat{C}$ is, as
is shown by the example of $\Gamma $-spaces
in~\cite{bousfield-friedlander}, where it is also shown that the
expected characterization of $\cat{S}$-local fibrations does not hold.

\section{The stable model structure}\label{sec-stable}

Our plan now is to apply Bousfield localization to the projective model
structure on $\BF $ to obtain a model structure with respect to which
$G$ is a Quillen equivalence.  In order to do this, we will have to
prove that the projective model structure makes $\BF $ into a cellular
model category when $\cat{C}$ is left proper cellular.  We will prove
this technical result in the appendix.  In this section, we will assume
that $\BF $ is cellular, find a good set $\cat{S}$ of maps to form the
stable model structure as the $\cat{S}$-localization of the projective model
structure, and prove that $G$ is a Quillen equivalence with respect to
the stable model structure. 

Just as in symmetric spectra~\cite{hovey-shipley-smith}, we want the
stable equivalences to be maps which induce isomorphisms on all
cohomology theories.  Cohomology theories will be represented by the
appropriate analogue of $\Omega $-spectra. 

\begin{definition}\label{defn-omega-spec}
A spectrum $X$ is a \emph{$U$-spectrum} if $X$ is level fibrant and
the adjoint $X_{n}\xrightarrow{\widetilde{\sigma }}UX_{n+1}$ of the
structure map is a weak equivalence for all $n$.  
\end{definition}

Of course, if $\cat{C}$ is the category of pointed simplicial sets or
pointed topological spaces, and $G$ is the suspension functor,
$U$-spectra are just $\Omega $-spectra.  We will find a set $\cat{S}$
of maps of $\BF $ such that the $\cat{S}$-local objects are the
$U$-spectra.  To do so, note that if $\map (A,X_{n})\xrightarrow{}\map
(A,UX_{n+1})$ is a weak equivalence of simplicial sets for all cofibrant
$A$ in $\cat{C}$, then $X_{n}\xrightarrow{}UX_{n+1}$ will be a weak
equivalence as required.  Since $\cat{C}$ is cofibrantly generated, we
should not need all cofibrant $A$, but only those $A$ related to the
generating cofibrations.  This is true, but the proof is somewhat
technical; the reader might be well-advised to skip the following
proof.  

\begin{proposition}\label{prop-weak-equivs-cofib-gen}
Suppose $\cat{C}$ is a left proper cofibrantly generated model category
with generating cofibrations $I$, and $f\mathcolon X\xrightarrow{}Y$ is
a map.  Then $f$ is a weak equivalence if and only if $\map
(C,X)\xrightarrow{}\map (C,Y)$ is a weak equivalence for all domains and
codomains $C$ of maps of $I$.  
\end{proposition}

\begin{proof}
The only if half is clear.  Since every cofibrant object is a retract of
a cell complex (\ie an object $A$ such that the map $0\xrightarrow{}A$
is a transfinite composition of pushouts of maps of $I$), it suffices to
show that $\map (A,f)$, or, equivalently, $\Map _{r}(A,Rf)$, is a weak
equivalence for all cell complexes $A$.  Given a cell complex $A$, there
is an ordinal $\lambda $ and a $\lambda $-sequence
\[
0=A_{0}\xrightarrow{} A_{1}\xrightarrow{} \dots \xrightarrow{} A_{\beta
} \xrightarrow{} \dots 
\]
with colimit $A_{\lambda }=A$.  We will show by transfinite induction on
$\beta $ that $\Map_{r} (A_{\beta },Rf)$ is a weak equivalence for all
$\beta \leq \lambda $.  Since $A_{0}=0$, getting started is easy.  For
the successor ordinal case of the induction, suppose $\map (A_{\beta
},f)$ is a weak equivalence.  We have a pushout square
\[
\begin{CD}
C @>>> A_{\beta } \\
@VgVV @VVi_{\beta }V \\
D @>>> A_{\beta +1}
\end{CD}
\]
where $g$ is a map of $I$.  Factor the composite
$QC\xrightarrow{}C\xrightarrow{}D$ into a cofibration
$\widetilde{g}\mathcolon QC\xrightarrow{}\widetilde{D}$ followed by a
trivial fibration $\widetilde{D}\xrightarrow{}D$.  In the terminology
of~\cite{hirschhorn}, $\widetilde{g}$ is a cofibrant approximation to
$g$.  By~\cite[Proposition~12.3.2]{hirschhorn}, there is a cofibrant
approximation $\widetilde{i_{\beta }}\mathcolon \widetilde{A_{\beta
}}\xrightarrow{}\widetilde{A_{\beta +1}}$ to $i_{\beta }$ which is a
pushout of $\widetilde{g}$.  This is where we need our model category to
be left proper.  For any fibrant object $Z$, the functor
$\Map _{r}(-,Z)$ converts colimits to limits and cofibrations to
fibrations~\cite[Corollary~5.4.4]{hovey-model}.  Hence we have two
pullback squares of fibrant simplicial sets
\[
\begin{CD}
\Map _{r}(\widetilde{A_{\beta +1}},RZ) @>>> \Map _{r}(\widetilde{D},RZ) \\
@VVV @VVV \\
\Map _{r}(\widetilde{A_{\beta }},RZ) @>>> \Map _{r}(QC,RZ)
\end{CD}
\]
where $Z=X$ and $Z=Y$, respectively.  Here the vertical maps are
fibrations.  There is a map from the square with $Z=X$ to the square
with $Z=Y$ induced by $f$.  By hypothesis, this map is a weak
equivalence on every corner except possibly the upper left.  But then
Dan Kan's cube lemma (see~\cite[Lemma 5.2.6]{hovey-model}, where the
dual of the version we need is proved, or~\cite{kan-model}) implies that
the map on the upper left corner is also a weak equivalence, and hence
that $\Map_{r} (A_{\beta +1},Rf)$ is a weak equivalence.

Now suppose $\beta $ is a limit ordinal and $\Map_{r} (A_{\gamma },Rf)$
is a weak equivalence for all $\gamma <\beta $.  Then, for $Z=X$ or
$Z=Y$, the simplicial sets $\Map _{r}(A_{\gamma },RZ)$ define a
limit-preserving functor $\beta^{\textup{op}} \xrightarrow{}\sset $ such
that each map $\Map _{r}(A_{\gamma +1},RZ)\xrightarrow{}\Map
_{r}(A_{\gamma },RZ)$ is a fibration of fibrant simplicial sets.  There
is a natural transformation from the functor with $Z=X$ to the functor
with $Z=Y$, and by hypothesis this map is a weak equivalence at every
stage.  It follows that it is a weak equivalence on the inverse limits,
as one can see in different ways.  The simplest is probably to note that
the inverse limit is a right Quillen
functor~\cite[Corollary~5.1.6]{hovey-model}.  Thus $\Map_{r} (A_{\beta
},Rf)$ is a weak equivalence, as required.  This completes the
transfinite induction and the proof.
\end{proof}

Note that the left properness assumption in
Proposition~\ref{prop-weak-equivs-cofib-gen} is unnecessary when the
domains of the generating cofibrations are themselves cofibrant, since
there is then no need to apply cofibrant approximation.  

In view of Proposition~\ref{prop-weak-equivs-cofib-gen}, we need to
choose our set $\cat{S}$ so as to make 
\[
\map (C,X_{n})\xrightarrow{}\map (C,UX_{n+1})
\]
a weak equivalence for all $\cat{S}$-local objects $X$ and all domains and
codomains $C$ of the generating cofibrations $I$.  Adjointness implies
that, if $X$ is level fibrant, $\map (C,X_{n})\cong \map (F_{n}QC,X)$ in
$\ho \sset $, since $F_{n}QC=(LF_{n})C$, where $LF_{n}$ is the
total left derived functor of $F_{n}$.  Also, $\map
(C,UX_{n+1})\cong \map (F_{n+1}GQC,X)$.  In view of this, we make the
following definition.

\begin{definition}\label{defn-S}
Suppose $\cat{C}$ is a left proper cellular model category with
generating cofibrations $I$, and $G$ is a left Quillen endofunctor of
$\cat{C}$.  Define the set $\cat{S}$ of maps in $\BF $ as $\{F_{n+1}GQC
\xrightarrow{s_{n}^{QC}} F_{n}QC \}$, as $C$ runs through the set of
domains and codomains of the maps of $I$ and $n$ runs through the
non-negative integers.  Here the map $s_{n}^{QC}$ is adjoint to the
identity map of $GQC$.  Define the \emph{stable model structure} on $\BF
$ to be the localization of the projective model structure on $\BF $
with respect to this set $\cat{S}$.  We refer to the $\cat{S}$-local
weak equivalences as \emph{stable equivalences}, and to the
$\cat{S}$-local fibrations as \emph{stable fibrations}.
\end{definition}

\begin{theorem}\label{thm-stable-equiv}
Suppose $\cat{C}$ is a left proper cellular model category and $G$ is a
left Quillen endofunctor of $\cat{C}$.  Then the stably fibrant objects
in $\BF$ are the $U$-spectra.  Furthermore, for all cofibrant $A\in
\cat{C}$ and for all $n\geq 0$, the map
$F_{n+1}GA\xrightarrow{s_{n}^{A}}F_{n}A$ is a stable equivalence.
\end{theorem}

\begin{proof}
By definition, $X$ is $\cat{S}$-local if and only if $X$ is level
fibrant and
\[
\map (F_{n}QC,X)\xrightarrow{}\map (F_{n+1}GQC,X)
\]
is a weak equivalence for all $n\geq 0$ and all domains and codomains
$C$ of maps of $I$.  By the comments preceding Definition~\ref{defn-S},
this is equivalent to requiring that $X$ be level fibrant and that the
map $\map (C,X_{n})\xrightarrow{}\map (C,UX_{n+1})$ be a weak
equivalence for all $n\geq 0$ and all domains and codomains $C$ of maps
of $I$.  By Proposition~\ref{prop-weak-equivs-cofib-gen}, this is
equivalent to requiring that $X$ be a $U$-spectrum.

Now, by definition, $s_{n}^{A}$ is a stable equivalence if and only if
$\map (s_{n}^{A},X)$ is a weak equivalence for all $U$-spectra $X$.  But
by adjointness, $\map (s_{n}^{A},X)$ can be identified with $\map
(A,X_{n})\xrightarrow{}\map (A,UX_{n+1})$.  Since
$X_{n}\xrightarrow{}UX_{n+1}$ is a weak equivalence between fibrant
objects, so is $\map (s_{n}^{A},X)$.  
\end{proof}

We would now like to claim that the stable model structure on $\BF $
that we have just defined is a generalization of the stable model
structure on spectra of topological spaces or simplicial sets defined
in~\cite{bousfield-friedlander}.  This cannot be a trivial observation,
however, both because our approach is totally different and because of
Remark~\ref{rem-subtle}.  

\begin{corollary}\label{cor-stable-spectra}
If $\cat{C}$ is either the category of pointed simplicial sets or
pointed topological spaces, and $G$ is the suspension functor given by
smashing with $S^{1}$, then the stable model structure on $\BF$
coincides with the stable model structure on the category of
Bousfield-Friedlander spectra~\cite{bousfield-friedlander}.
\end{corollary}

\begin{proof}
We know already that the cofibrations are the same in the stable model
structure on $\BF $ and the stable model structure
of~\cite{bousfield-friedlander}.  We will show that the weak
equivalences are the same.  In any model category at all, a map $f$ is a
weak equivalence if and only if $\map (f,X)$ is a weak equivalence of
simplicial sets for all fibrant $X$.  Construction of $\map (f,X)$
requires replacing $f$ by a cofibrant approximation $f'$ and building
cosimplicial resolutions of the domain and codomain of $f'$.  In the
case at hand, we can do the cofibrant replacement and build the
cosimplicial resolutions in the strict model category of spectra, since
the cofibrations do not change under localization.  Thus $\map (f,X)$ is
the same in both the stable model structure on $\BF$ and in the stable
model category of Bousfield and Friedlander.  Since the stably fibrant
objects are also the same, the corollary holds.
\end{proof}

The purpose of the stable model structure is to make the prolongation of
$G$ into a Quillen equivalence.  We begin the process of proving this
with the following corollary. 

\begin{corollary}\label{cor-G-Quillen}
Suppose $\cat{C}$ is a left proper cellular model category and $G$ is a
left Quillen endofunctor of $\cat{C}$.  Then the prolongation of $G$ to
a functor $G\mathcolon \BF\xrightarrow{}\BF$ is a Quillen functor with
respect to the stable model structure.  
\end{corollary}

\begin{proof}
In view of Hirschhorn's localization theorem~\ref{thm-localization}, we
must show that $G(Qf)$ is a stable equivalence for all $f\in \cat{S}$.
Since the domains and codomains of the maps of $\cat{S}$ are already
cofibrant, it is equivalent to show that $Gf$ is a stable equivalence
for all $f\in S$.  Since $GF_{n}=F_{n}G$, we have
$Gs_{n}^{A}=s_{n}^{GA}$.  In view of Theorem~\ref{thm-stable-equiv},
this map is a weak equivalence whenever $A$, and hence $GA$, is
cofibrant.  Taking $A=QC$, where $C$ is a domain or codomain of a map of
$I$, completes the proof.
\end{proof}

We will now show that $G$ is in fact a Quillen \emph{equivalence} with
respect to the stable model structure.  To do so, we introduce the shift
functors.

\begin{definition}\label{defn-shift}
Suppose $\cat{C}$ is a model category and $G$ is a left Quillen
endofunctor of $\cat{C}$.  Define the \emph{shift functors} $t\mathcolon
\BF \xrightarrow{}\BF $ and $s\mathcolon
\BF \xrightarrow{}\BF $ by $(sX)_{n}=X_{n+1}$ and
$(tX)_{n}=X_{n-1}$, $(tX)_{0}=0$, with the evident structure maps.  Note
that $t$ is left adjoint to $s$.  
\end{definition}

It is clear that $s$ preserves level equivalences and level fibrations,
so $t$ is a left Quillen functor with respect to the strict model
structure on $\BF $, and $s$ is a right
Quillen functor.  Also, $sU=Us$, so $tG=Gt$.  Similarly, $\Ev _{n}s=\Ev
_{n+1}$, so $tF_{n}=F_{n+1}$.  It follows that $ts_{n}^{A}=s_{n+1}^{A}$,
so that $t$ is a Quillen functor with respect to the stable model
structure as well.  

We have now come to the main advantage of our approach to spectra.  We
can test whether a spectrum $X$ is a $U$-spectrum by checking that $X$
is level fibrant and by checking that the map $X\xrightarrow{}sUX$,
adjoint to the structure map of $X$, is a level equivalence.  The
analogous statement is \textbf{false} in the Bousfield-Friedlander
category~\cite{bousfield-friedlander}, because the extra twist map they
use (see Remark~\ref{rem-subtle}) means that there is \textbf{no map of
spectra} $X\xrightarrow{}sUX$ adjoint to the structure map of $X$!  Our
interpretation of this is that the Bousfield-Friedlander approach, while
excellent at what it does, is probably not the right general
construction.  Further evidence for this is provided by the extreme
simplicity of the following proof. 

\begin{theorem}\label{thm-G-Quillen-equiv}
Suppose $\cat{C}$ is a left proper cellular model category and $G$ is a
left Quillen endofunctor of $\cat{C}$.  Then the functors $G\mathcolon
\BF \xrightarrow{}\BF $ and $t\mathcolon \BF \xrightarrow{}\BF $ are
Quillen equivalences with respect to the stable model structures.
Furthermore, $Rs$ is naturally isomorphic to $LG$, and $RU$ is
naturally isomorphic to $Lt$.
\end{theorem}

\begin{proof}
There is a a natural map $X\xrightarrow{}sUX$ which is a weak
equivalence when $X$ is a stably fibrant object of $\BF $.  This means
that the total right derived functor $R(sU)$ is naturally isomorphic to
the identity functor on $\ho \BF $ (where we use the stable model
structure).  On the other hand, $R(sU)$ is naturally isomorphic to
$Rs\circ RU$ and also to $RU\circ Rs$, since $s$ and $U$ commute with
each other.  Thus the natural isomorphism from the identity to $R(sU)$
gives rise to an natural isomorphism $1\xrightarrow{}Rs\circ RU$ and a
natural isomorphism $RU\circ Rs\xrightarrow{}1$, displaying $RU$ and
$Rs$ as adjoint equivalences of categories.  It follows that $U$ and
$s$ are both Quillen equivalences, as required.  Since $RU$ is adjoint
to $LG$ and $Rs$ is adjoint to $Lt$, we must also have $RU$ naturally
isomorphic to $Lt$ and $Rs$ naturally isomorphic to $LG$.
\end{proof}

\section{The almost finitely generated case}\label{sec-fin-gen}

The reader may well object at this point that we have defined the stable
model structure on $\BF $ without ever defining stable homotopy groups.
This is because stable homotopy groups do not detect stable equivalences
in general.  The usual simplicial and topological situation is very
special.  The goal of this section is to put some hypotheses on
$\cat{C}$ and $G$ so that the stable model structure on $\BF $ behaves
similarly to the stable model structure on ordinary simplicial spectra.
In particular, we show that, if $\cat{C}$ is almost finitely generated
(defined below), the usual $Q$ construction gives a stable fibrant
replacement functor; thus, a map $f$ is a stable equivalence if and only
if $Qf$ is a level equivalence.  This allows us to characterize $\ho \BF
(F_{0}A,X)$ for well-behaved $A$ as the usual sort of colimit $\colim
\ho \cat{C}(G^{n}A,X_{n})$.  It also allows us to prove that the stable
model structure is right proper, so we get the expected characterization
of stable fibrations.

Most of the results in this section do not depend on the existence of
the stable model structure on $\BF $, so we do not usually need to assume
$\cat{C}$ is left proper cellular.  

We now define almost finitely generated model categories, as suggested
to the author by Voevodsky.  

\begin{definition}\label{def-finite}
An object $A$ of a category $\cat{C}$ is called \emph{finitely
presented} if the functor $\cat{C}(A,-)$ preserves direct limits of
sequences $X_{0}\xrightarrow{}X_{1}\xrightarrow{}\dots
\xrightarrow{}X_{n}\xrightarrow{}\dots $.  A cofibrantly generated model
category $\cat{C}$ is said to be \emph{finitely generated} if the
domains and codomains of the generating cofibrations and the generating
trivial cofibrations are finitely presented.  A cofibrantly generated
model category is said to be \emph{almost finitely generated} if the
domains and codomains of the generating cofibrations are finitely
presented, and if there is a set of trivial cofibrations $J'$ with
finitely presented domains and codomains such that a map $f$ \emph{whose
codomain is fibrant} is a fibration if and only if $f$ has the \rlp
$J'$.  
\end{definition}

This definition differs slightly from other definitions.  In particular,
an object $A$ is usually said to be finitely presented if $\cat{C}(A,-)$
preserves all directed (or, equivalently, filtered) colimits.  We are
trying to assume the minimum necessary.  Finitely generated model
categories were introduced in~\cite[Section~7.4]{hovey-model}, but in
that definition we assumed only that $\cat{C}(A,-)$ preserves
(transfinitely long) direct limits of sequences of \emph{cofibrations}.
The author would now prefer to call such model categories
\emph{compactly generated}.  Thus, the model category of simplicial sets
is finitely generated, but the model category on topological spaces is
only compactly generated.  Since we will only be working with (almost)
finitely generated model categories in this section, our results will
not apply to topological spaces.  We will indicate where
our results fail for compactly generated model categories, and a
possible way to amend them in the compactly generated case.

The definition of an almost finitely generated model category was
suggested by Voevodsky.  The problem with finitely generated, or,
indeed, compactly generated, model categories is that they are not
preserved by localization.  That is, if $\cat{C}$ is a finitely
generated left proper cellular model category, and $S$ is a set of maps,
then the Bousfield localization $L_{S}\cat{C}$ will not be finitely
generated, because we lose all control over the generating trivial
cofibrations in $L_{S}\cat{C}$.  However, if $S$ is a set of
cofibrations such that $X\otimes K$ is finitely presented for every
domain or codomain $X$ of a map of $S$ and every finite simplicial set
$K$, then $L_{S}\cat{C}$ will still be almost finitely generated.  (We
use a framing on $\cat{C}$ to construct $X\otimes K$).  Indeed, the
horns
\[
(A\otimes \Delta [n])\amalg _{A\otimes \Lambda ^{k}[n]}(B\otimes \Lambda
^{k}[n])\xrightarrow{}B\otimes \Delta [n]
\]
on the maps $A\xrightarrow{}B$ of $S$ are used to detect $S$-local
fibrant objects, and an $S$-local fibration between $S$-local fibrant
objects is just an ordinary fibration.  We can therefore take the set
$J'$ to consist of the horns on the maps of $S$ together with the old
set of generating trivial cofibrations.

In particular, Voevodsky has informed the author that he can make an
unstable motivic model category that is almost finitely generated, using
this approach.  For the reader's benefit, we summarize his construction.
The category $\cat{C}$ is the category of simplicial presheaves (of
sets) on the category of smooth schemes over some base scheme $k$.
There is a projective model structure on this category, where weak
equivalences and fibrations are defined objectwise from weak
equivalences and fibrations of simplicial sets.  The projective model
structure is finitely generated (using the fact that smooth schemes over
$k$ is an essentially small category).  There is an embedding of smooth
schemes into $\cat{C}$ as representable functors.  We need to localize
this model structure to take into account both the Nisnevich topology
and the fact that the functor $X\mapsto X\times \mathbb{A}^{1}$ should
be a Quillen equivalence.  To do so, we define a set $S'$ to consist of
the maps $X\times \mathbb{A}^{1}\xrightarrow{}X$ for every smooth scheme
$X$ and maps $P\xrightarrow{}X$ for every pullback square of smooth
schemes 
\[
\begin{CD}
B @>>> Y \\
@VVV @VVpV \\
A @>>j> X
\end{CD}
\] 
where $p$ is etale, $j$ is an open embedding, and
$p^{-1}(X-A)\xrightarrow{}X-A$ is an
isomorphism.  Here $P$ is the mapping cylinder $(B\amalg Y)\amalg
_{B\amalg B}(A\times \Delta [1])$.  We then define $S$ to consist of
mapping cylinders on the maps of $S'$.  The maps of $S$ are then
cofibrations whose domains and codomains are finitely presented (and
remain so after tensoring with any finite simplicial set), so the
Bousfield localization will be almost finitely generated.

There is then some work involving properties of the Nisnevich topology
to show that this model category is equivalent to the Morel-Voevodsky
motivic model category of~\cite{morel-voevodsky}, and to the model
category used by Jardine~\cite{jardine}.  

The reason we need almost finitely generated model categories is
because, in an almost finitely generated model category $\cat{C}$,
sequential colimits preserve trivial fibrations, fibrant objects, and
fibrations between fibrant objects.  Indeed, suppose we have a map of
sequences $p_{n}\mathcolon
X_{n}\xrightarrow{}Y_{n}$ that is a fibration between fibrant objects
for all $n$.  We show $\colim X_{n}$ is fibrant by testing that 
$X\xrightarrow{}0$ has the \rlp
$J'$.  We then test that $\colim p_{n}$ is a fibration by testing that
it has the \rlp $J'$.  The proof that sequential colimits preserve
trivial fibrations is similar.  

Now, given a spectrum $X$, there is an obvious candidate for a stable
fibrant replacement of $X$.  

\begin{definition}\label{defn-Q}
Suppose $G$ is a left Quillen endofunctor of a model category $\cat{C}$
with right adjoint $U$.  Define $R\mathcolon \BF \xrightarrow{}\BF $ to
be the functor $sU$, where $s$ is the shift functor.  Then we have a
natural map $\iota _{X}\mathcolon X\xrightarrow{}RX$, and we define
\[
R^{\infty }X = \colim (X\xrightarrow{\iota _{X}}RX\xrightarrow{R\iota
_{X}}R^{2}X\xrightarrow{R^{2}\iota _{X}}\dots \xrightarrow{R^{n-1}\iota
_{X}}R^{n}X\xrightarrow{R^{n}\iota _{X}}\dots ) .
\]
Let $j_{X}\mathcolon X\xrightarrow{}R^{\infty }X$ denote the obvious
natural transformation.  
\end{definition}

The following lemma, though elementary, is crucial.  

\begin{lemma}\label{lem-twist}
The maps $\iota _{RX}, R\iota _{X}\mathcolon RX\xrightarrow{}R^{2}X$
coincide.  
\end{lemma}

\begin{proof}
The map $\iota _{RX}$ is the adjoint of the structure map 
\[
GUX_{n+1}\xrightarrow{\varepsilon } X_{n+1}\xrightarrow{\eta } UGX_{n+1}
\xrightarrow{U\sigma } UX_{n+2}
\]
of $RX$, where $\varepsilon $ denotes the counit of the adjunction,
$\eta $ denotes the unit, and $\sigma $ denotes the structure map of
$X$.  Thus $\iota _{RX}$ is the composite 
\[
UX_{n+1} \xrightarrow{\eta } UGUX_{n+1} \xrightarrow{U\varepsilon }
UX_{n+1} \xrightarrow{U\eta } U^{2}GX_{n+1} \xrightarrow{U^{2}\sigma }
U^{2}X_{n+2}.  
\]
Since $U\varepsilon \circ \eta $ is the identity, it follows that $\iota
_{RX}=R\iota _{X}$.  
\end{proof}

We stress that Lemma~\ref{lem-twist} fails for symmetric spectra, and it
is the major reason we must work with finitely generated model
categories rather than compactly generated model categories.  Indeed,
in the compactly generated case, $R^{\infty }$ is not a good functor,
since maps out of one of the domains of the generating cofibrations will
not preserve the colimit that defines $R^{\infty }$.  The obvious thing
to try is to replace the functor $R$ by a functor $W$, obtained by
factoring $X\xrightarrow{}RX$ into a projective cofibration
$X\xrightarrow{}WX$ followed by a level trivial fibration
$WX\xrightarrow{}RX$.  The difficulty with this plan is that we do not
see how to prove Lemma~\ref{lem-twist} for $W$.  An alternative plan
would be to use the mapping cylinder $X\xrightarrow{}W'X$ on
$X\xrightarrow{}RX$; this might make Lemma~\ref{lem-twist} easier to
prove, but the map $X\xrightarrow{}W'X$ will not be a cofibration.  The
map $X\xrightarrow{}W'X$ may, however, be good enough for the required
smallness properties to hold.  It is a closed inclusion if
$\cat{C}$ is topological spaces, for example. The author knows of no
good general theorem in the compactly generated case.  

This lemma leads immediately to the following proposition.

\begin{proposition}\label{prop-strict}
Suppose $G$ is a left Quillen endofunctor of a model category $\cat{C}$,
and suppose that its right adjoint $U$ preserves sequential colimits.
Then the map $\iota _{R^{\infty }X}\mathcolon R^{\infty
}X\xrightarrow{}R(R^{\infty }X)$ is an isomorphism.  In particular, if
$X$ is level fibrant, $\cat{C}$ is almost finitely generated, and $U$
preserves sequential colimits, then $R^{\infty }X$ is a $U$-spectrum.
\end{proposition}

\begin{proof}
The map $\iota _{R^{\infty }X}$ is the colimit of the vertical maps in
the diagram below. 
\[
\begin{CD}
X @>\iota _{X}>> RX @>R\iota _{X}>> R^{2}X @>R^{2}\iota _{X}>> \dots
@>R^{n-1}\iota _{X}>> R^{n}X @>R^{n}\iota _{X}>> \dots \\
@V\iota _{X}VV @V\iota _{RX}VV @V\iota _{R^{2}X}VV @. @V\iota
_{R^{n}X}VV \\
RX @>>R\iota _{X}> R^{2}X @>>R^{2}\iota _{X}> R^{3}X @>>R^{3}\iota _{X}>
\dots @>>R^{n}\iota _{X}> R^{n+1}X @>>R^{n+1}\iota _{X}> \dots 
\end{CD}
\]
Since the vertical and horizontal maps coincide, the result follows.
For the second statement, we note that if $X$ is level fibrant, each
$R^{n}X$ is level fibrant since $R$ is a right Quillen functor (with
respect to the projective model structure).  Since sequential colimits
in $\cat{C}$ preserve fibrant objects, $R^{\infty }X$ is level fibrant,
and hence a $U$-spectrum.
\end{proof}

\begin{proposition}\label{prop-R-is-level}
Suppose $G$ is a left Quillen endofunctor of a model category $\cat{C}$
with right adjoint $U$.  If $\cat{C}$ is almost finitely generated, and
$X$ is a $U$-spectrum, then the map $j_{X}\mathcolon
X\xrightarrow{}R^{\infty }X$ is a level equivalence.
\end{proposition}

\begin{proof}
By assumption, the map $\iota _{X}\mathcolon X\xrightarrow{}RX$ is a
level equivalence between level fibrant objects.  Since $R$ is a right
Quillen functor, $R^{n}\iota _{X}$ is a level equivalence as well.  Then
the method of~\cite[Corollary~7.4.2]{hovey-model} completes the proof.
Recall that this method is to use factorization to construct a sequence
of projective trivial cofibrations $Y_{n}\xrightarrow{}Y_{n+1}$ with
$Y_{0}=X$ and a level trivial fibration of sequences
$Y_{n}\xrightarrow{}R^{n}X$.  Then the map $X\xrightarrow{}\colim Y_{n}$
will be a projective trivial cofibration.  Since sequential colimits in
$\cat{C}$ preserve trivial fibrations, the map $\colim
Y_{n}\xrightarrow{}R^{\infty }X$ will still be a level trivial
fibration.
\end{proof}

Proposition~\ref{prop-R-is-level} gives us a slightly better method of
detecting stable equivalences.  

\begin{corollary}\label{cor-detecting}
Suppose $G$ is a left Quillen endofunctor of a model category $\cat{C}$
with right adjoint $U$.  Suppose $\cat{C}$ is almost finitely generated
and $U$ preserves sequential colimits.  Then a map $f\mathcolon
A\xrightarrow{}B$ is a stable equivalence in $\BF $ if and only if $\map
(f,X)$ is a weak equivalence for all level fibrant spectra $X$ such that
$\iota _{X}\mathcolon X\xrightarrow{}RX$ is an isomorphism.
\end{corollary}

\begin{proof}
By definition, $f$ is a stable equivalence if and only if $\map (f,Y)$
is a weak equivalence for all $U$-spectra $Y$.  But we have a level
equivalence $Y\xrightarrow{}R^{\infty }Y$ by
Proposition~\ref{prop-R-is-level}, and so it suffices to know 
that $\map (f, R^{\infty }Y)$ is a weak equivalence for all $U$-spectra
$Y$.  But, by Proposition~\ref{prop-strict}, $\iota _{R^{\infty }Y}$ is
an isomorphism.  
\end{proof}

This corollary, in turn, allows us to prove that $R^{\infty }$ detects
stable equivalences. The following theorem is similar
to~\cite[Theorem~3.1.11]{hovey-shipley-smith}.  

\begin{theorem}\label{thm-R-detects}
Suppose $G$ is a left Quillen endofunctor of a model category $\cat{C}$
with right adjoint $U$.  Suppose that $\cat{C}$ is almost finitely
generated and sequential colimits in $\cat{C}$ preserve finite products.
Suppose also that $U$ preserves sequential colimits.  If $f\mathcolon
A\xrightarrow{}B$ is a map in $\BF $ such that $R^{\infty }f$ is a level
equivalence, then $f$ is a stable equivalence.
\end{theorem}

\begin{proof}
Suppose $X$ is a $U$-spectrum such that the map $\iota _{X}\mathcolon
X\xrightarrow{}RX$ is an isomorphism.  We will show that $\map (f,X)$ as
a retract of $\map (R^{\infty }f,R^{\infty }X)$; this will obviously
complete the proof.  We first note that there is a natural map $\map
(R^{\infty }C,R^{\infty }X)\xrightarrow{}\map (C,X)$ obtained by
precomposition with $C\xrightarrow{}R^{\infty }C$ and postcomposition
with $k\mathcolon R^{\infty }X\xrightarrow{}X$.  Here $k$ is the inverse
of the map $j_{X}\mathcolon X\xrightarrow{}R^{\infty }X$, which is an
isomorphism since $\iota _{X}$ is so.  On the other hand, we claim that
there is also a natural map $\map (C,X)\xrightarrow{}\map (R^{\infty
}C,R^{\infty }X)$ obtained by applying the total right derived functor
of $R^{\infty }$.  This is not obvious, since we are asserting that the
total right derived functor of $R^{\infty }$ preserves the enrichment of the
projective homotopy category of $\BF $ over the homotopy category of
simplicial sets, even though $R^{\infty }$ is not a right Quillen
functor.  Nevertheless, if we assume that this natural map exists, it
follows easily that the composite $\map (C,X)\xrightarrow{}\map
(R^{\infty }C,R^{\infty }X)\xrightarrow{}\map (C,X)$ is the identity (in
the homotopy category of simplicial sets), and therefore that $\map
(f,X)$ is a retract of $\map (R^{\infty }f,R^{\infty }X)$.

It remains to show that the total right derived functor of $R^{\infty }$
preserves the enriched structure.  We first point out that $R^{\infty }$
preserves level fibrations between level fibrant objects and all level
trivial fibrations, because $\cat{C}$ is almost finitely generated.  It
follows from Ken Brown's lemma~\cite[Lemma~1.1.12]{hovey-model} that
$R^{\infty }$ preserves level equivalences between level fibrant
objects.  Because sequential colimits in $\cat{C}$ preserve finite
products, $R^{\infty }$ also preserves finite products.

We claim that, for any functor $H$ on a model category $\cat{D}$ that
preserves fibrations between fibrant objects, weak equivalences between
fibrant objects, and finite products, the total right derived functor of
$H$ preserves the enriched structure over $\ho \sset $.  Indeed,
analysis of the definition of this
enrichment~\cite[Chapter~5]{hovey-model} shows that it suffices to check
that such a functor $H$ preserves simplicial frames on a fibrant object
$Y$.  A simplicial frame on $Y$ is a factorization $\ell _{\bullet
}Y\xrightarrow{\alpha }Y_{*}\xrightarrow{\beta }r_{\bullet }Y$ in the
diagram category of simplicial objects in $\cat{D}$, where $\ell
_{\bullet }Y$ is the constant simplicial diagram on $Y$, the $n$th space
of $r_{\bullet }Y$ is the product $Y^{n+1}$, $\alpha $ is a level
equivalence, and $\beta $ is a level fibration.  We further require that
$\beta $ is an isomorphism in degree $0$.  Since $H$ preserves products,
weak equivalences between fibrant objects, and fibrations between
fibrant objects, it follows that $HY_{*}$ is a simplicial frame on $HY$.
\end{proof}

\begin{corollary}\label{cor-j-is-stable-equiv}
Suppose $G$ is a left Quillen endofunctor of a model category $\cat{C}$
with right adjoint $U$.  Suppose that $\cat{C}$ is almost finitely
generated and sequential colimits in $\cat{C}$ preserve finite products.
Suppose also that $U$ preserves sequential colimits.  Then
$j_{A}\mathcolon A\xrightarrow{}R^{\infty }A$ is a stable equivalence
for all $A\in \BF $.
\end{corollary}

\begin{proof}
One can easily check that $R^{\infty }j_{A}$ is an isomorphism, using
Proposition~\ref{prop-strict}.  
\end{proof}

Finally, we get the desired characterization of stable equivalences. 

\begin{theorem}\label{thm-stable}
Suppose $G$ is a left Quillen endofunctor of a model category $\cat{C}$
with right adjoint $U$.  Suppose that $\cat{C}$ is almost finitely
generated and sequential colimits in $\cat{C}$ preserve finite products.
Suppose as well that $U$ preserves sequential colimits.  Let $L'$ denote
a fibrant replacement functor in the projective model structure on $\BF
$.  Then, for all $A\in \BF $, the map $A\xrightarrow{}R^{\infty }L'A$
is a stable equivalence into a $U$-spectrum.  Also, a map $f\mathcolon
A\xrightarrow{}B$ is a stable equivalence if and only if $R^{\infty
}L'f$ is a level equivalence.
\end{theorem}

\begin{proof}
The first statement follows immediately from
Proposition~\ref{prop-strict} and
Corollary~\ref{cor-j-is-stable-equiv}.  By the first statement, if $f$
is a stable equivalence, so is $R^{\infty }L'f$.  Since $R^{\infty }L'f$
is a map between $U$-spectra, it is a stable equivalence if and only if
it is a level equivalence.  The converse follows from
Theorem~\ref{thm-R-detects}. 
\end{proof}

Since we did not need the existence of the stable model structure to
prove Theorem~\ref{thm-stable}, one can imagine attempting to construct
it from the functor $R^{\infty }L'$.  This is, of course, the original
approach of Bousfield-Friedlander~\cite{bousfield-friedlander}, and this
approach has been generalized by Schwede~\cite{schwede}.  Also, if one
has some way to detect level equivalences in $\cat{C}$, say using
appropriate generalizations of homotopy groups, Theorem~\ref{thm-stable}
implies that stable equivalences in $\BF $ are detected by the
appropriate generalizations of stable homotopy groups.  One can see
these generalizations in the following corollary as well.  

\begin{corollary}\label{cor-maps-stable}
Suppose $\cat{C}$ is a pointed, left proper, cellular, almost finitely
generated model category where sequential colimits preserve finite
products.  Suppose $G\mathcolon \cat{C}\xrightarrow{}\cat{C}$ is a left
Quillen functor whose right adjoint $U$ commutes with sequential
colimits.  Finally, suppose $A$ is a finitely presented cofibrant object
of $\cat{C}$ that has a finitely presented cylinder object $A\times I$.
Then
\[
\ho \BF (F_{k}A,Y) =\colim_{m} \ho \cat{C}(A,U^{m}Y_{k+m}).  
\]
for all level fibrant $Y\in \BF $.  
\end{corollary}

Here we are using the stable model structure to form $\ho \BF $, of
course.  

\begin{proof}
We have $\ho \BF (F_{k}A, Y)=\BF (F_{k}A,R^{\infty }Y)/\sim $, by
Theorem~\ref{thm-stable}, where $\sim $ denotes the left homotopy
relation.  We can use the cylinder object $F_{k}(A\times I)$ as the
source for our left homotopies.  Then adjointness implies that $\BF
(F_{k}A, R^{\infty }Y)/\sim =\cat{C}(A, \Ev _{k}R^{\infty }Y)/\sim $.
Since $A$ and $A\times I$ are finitely presented, we get the required
result.  
\end{proof}

By assuming slightly more about $\cat{C}$, we can also characterize the
stable fibrations.  

\begin{corollary}\label{cor-right-proper}
Suppose $\cat{C}$ is a pointed, proper, cellular, almost finitely
generated model category such that sequential colimits preserve
pullbacks.  Suppose $G\mathcolon \cat{C}\xrightarrow{}\cat{C}$ is a left
Quillen functor whose right adjoint $U$ commutes with sequential
colimits.  Then the stable model structure on $\BF $ is proper.  In
particular, a map $f\mathcolon X\xrightarrow{}Y$ is a stable fibration
if and only if $f$ is a level fibration and the diagram
\[
\begin{CD}
X @>>> R^{\infty }L'X \\
@VfVV @VVLfV \\
Y @>>> R^{\infty }L'Y
\end{CD}
\]
is a homotopy pullback square in the projective model structure, where
$L'$ is a fibrant replacement functor in the projective model
structure. 
\end{corollary}

\begin{proof}
We wil actually show that, if $p\mathcolon X\xrightarrow{}Y$ is a level
fibration and $f\mathcolon B\xrightarrow{}Y$ is a stable equivalence,
the pullback $B\times _{A}Y\xrightarrow{}X$ is a stable equivalence.
The first step is to use the right properness of the projective model
structure on $\BF $ to reduce to the case where $B$ and $Y$ are level
fibrant.  Indeed, let $Y'=L'Y$, $B'=L'B$, and $f'=L'f$.  Then factor the
composite $X\xrightarrow{}Y\xrightarrow{}Y'$ into a projective trivial
cofibration $X'\xrightarrow{}Y'$ followed by a level fibration
$p'\mathcolon X'\xrightarrow{}Y'$.  Then we have the commutative diagram
below, 
\[
\begin{CD}
B @>f>> Y @<p<< X \\
@VVV @VVV @VVV \\
B' @>>f'> Y' @<<p''< X'
\end{CD}
\]
where the vertical maps are level equivalences.  Then Proposition~12.2.4
and Corollary~12.2.8 of~\cite{hirschhorn}, which depend on the
projective model structure being right proper, imply that the induced
map $B\times _{Y}X\xrightarrow{}B'\times _{Y'}X'$ is a level
equivalence.  Hence $B\times _{Y}X\xrightarrow{}X$ is a stable
equivalence if and only if $B'\times _{Y'}X'\xrightarrow{}X'$ is a
stable equivalence, and so we can assume $B$ and $X$ are level fibrant.  

Now let $S$ denote the pullback square below.  
\[
\begin{CD}
B\times _{Y}X @>>> X \\
@VVV @VVpV \\
B @>>f> Y
\end{CD}
\]
Then $R^{n}S$ is a pullback square for all $n$, and there are maps
$R^{n}S\xrightarrow{R^{n}\iota _{S}}R^{n+1}S$.  Since pullbacks commute
with sequential colimits, $R^{\infty }S$ is a pullback square.
Furthermore, $R^{\infty }p$ is a level fibration, since sequential
colimits in $\cat{C}$ preserve fibrations between level fibrant objects.
Since $f$ is a stable equivalence between level fibrant spectra,
$R^{\infty }f$ is a level equivalence by Theorem~\ref{thm-stable}.  So,
since the projective model structure is right proper, the map $R^{\infty
}(B\times _{Y}X\xrightarrow{}X)$ is a level equivalence, and thus
$B\times _{Y}X\xrightarrow{}X$ is a stable equivalence.

The characterization of stable fibrations then follows
from~\cite[Proposition~3.6.8]{hirschhorn}.  
\end{proof}

\section{Properties of the stabilization functor}\label{sec-funct}

In this section we explore some of the properties of the correspondence
$(\cat{C},G) \mapsto \BF $, where, throughout this section, we mean the
stable model structure on $\BF $.  We begin by showing that, if $G$ is
already a Quillen equivalence, then the embedding
$\cat{C}\xrightarrow{}\BF $ is a Quillen equivalence.  This important
fact is as close as we can get to proving that $\BF $ is the initial, up
to homotopy, stabilization of $\cat{C}$ with respect to $G$.  We also
show that $\BF $ is functorial in the pair $(\cat{C},G)$, with a
suitable definition of maps of pairs.  Under mild hypotheses, we show
that $\BF $ preserves Quillen equivalences in the pair $(\cat{C},G)$.
In particular, this means, for example, that the Bousfield-Friedlander
category of spectra of simplicial sets does not depend, up to Quillen
equivalence, on which model of the circle $S^{1}$ one chooses.  We
conclude the section by pointing out that our stabilization condition
preserves some monoidal structure.  For example, if $\cat{C}$ is a
simplicial model category, and $G$ is a simplicial functor, then $\BF $
is a gain a simplicial model category, and the extension of $G$ is again
a simplicial functor.  However, if $\cat{C}$ is monoidal, and $G$ is a
monoidal functor, $\BF $ will almost never be a monoidal category; this
is the reason we need the symmetric spectra introduced in the next
section.  

\begin{theorem}\label{thm-already-equiv}
Suppose $\cat{C}$ is a left proper cellular model category, and suppose
$G$ is a left Quillen endofunctor of $\cat{C}$ that is a Quillen
equivalence.  Then $F_{0}\mathcolon \cat{C}\xrightarrow{}\BF $
is a Quillen equivalence, where $\BF $ has the stable model structure. 
\end{theorem}

\begin{proof}
We first point out that the right adjoint $\Ev _{0}\mathcolon \BF
\xrightarrow{}\cat{C}$ reflects weak equivalences between fibrant
objects.  Indeed, suppose $X$ and $Y$ are $U$-spectra, and $f\mathcolon
X\xrightarrow{}Y$ is a map such that $\Ev _{0}f=f_{0}$ is a weak
equivalence.  Then, because $X$ and $Y$ are $U$-spectra, $U^{n}f_{n}$ is
a weak equivalence for all $n$.  Since $G$ is a Quillen equivalence, $U$
reflects weak equivalences between fibrant objects
by~\cite[Corollary~1.3.16]{hovey-model}.  Thus $f_{n}$ is a weak
equivalence for all $n$, and so $f$ is a level equivalence and hence a
stable equivalence, as required.

In view of~\cite[Corollary~1.3.16]{hovey-model}, to complete the proof
it suffices to show that $X\xrightarrow{}(L_{\cat{S}}F_{0}X)_{0}$ is a
weak equivalence for all cofibrant $X\in \cat{C}$, where $L_{\cat{S}}$
denotes a stably fibrant replacement functor in $\BF $.  Let $R'$ denote
a fibrant replacement functor in the strict model structure on $\BF $.
Then $X\xrightarrow{}(R'F_{0}X)_{0}$ is certainly a weak equivalence.
We claim that $R'F_{0}X$ is already a $U$-spectrum.  Suppose for the
moment that this is true; then by lifting we can construct a stable
equivalence $R'F_{0}X\xrightarrow{}L_{\cat{S}}F_{0}X$, and since
$R'F_{0}X$ is a $U$-spectrum, this map is in fact a level equivalence.
Hence the map $X\xrightarrow{}(L_{\cat{S}}F_{0}X)_{0}$ is a weak
equivalence, as required.

It remains to prove that $R'F_{0}X$ is a $U$-spectrum.  Since $G$ is a
Quillen equivalence, the map
$(F_{0}X)_{n}=G^{n}X\xrightarrow{}URG^{n+1}X=UR(F_{0}X)_{n+1}$ is a weak
equivalence.  By lifting, we can factor the weak equivalence
$(F_{0}X)_{n+1}\xrightarrow{}(R'F_{0}X)_{n+1}$ through the trivial
cofibration $(F_{0}X)_{n+1}\xrightarrow{}R(F_{0}X)_{n+1}$.  This implies
that the map $(F_{0}X)_{n}\xrightarrow{}U(R'F_{0}X)_{n+1}$ is a weak
equivalence, and hence that $R'F_{0}X$ is a $U$-spectrum, as required.
\end{proof}

In particular, this theorem means that the passage $(\cat{C},G)\mapsto
(\BF ,G)$ is idempotent, up to Quillen equivalence.  This suggests that
we are doing some kind of fibrant replacement of $(\cat{C},G)$, but the
author knows no way of making this precise. 

We now examine the functoriality of the stable model structure on $\BF $.  

\begin{definition}\label{defn-map-of-pairs}
Suppose $\cat{C}$ and $\cat{D}$ are left proper cellular model
categories, $G$ is a left Quillen endofunctor of $\cat{C}$, and $H$ is a
left Quillen endofunctor of $\cat{D}$.  A \emph{map of pairs} $(\Phi
,\tau ) \mathcolon (\cat{C},G) \xrightarrow{} (\cat{D},H)$ is a left
Quillen functor $\Phi \mathcolon \cat{C}\xrightarrow{} \cat{D}$ and a
natural transformation $\tau \mathcolon \Phi G \xrightarrow{} H\Phi $
such that $\tau _{A}$ is a weak equivalence for all cofibrant $A\in
\cat{C}$.
\end{definition}

Note that there is an obvious associative and unital composition of maps
of pairs.  

\begin{proposition}\label{prop-functor}
Suppose $(\Phi ,\tau )\mathcolon (\cat{C},G) \xrightarrow{} (\cat{D},H)$
is a map of pairs.  Then there is an induced map of pairs $(\funcBF{\Phi
}, \funcBF{\tau }) \mathcolon (\BF ,G)
\xrightarrow{}(\genBF{\cat{D}}{H},H)$ such that $\funcBF{\Phi } \circ
F_{n}=F_{n}\Phi $.  This induced map of pairs is compatible with
composition and identities.
\end{proposition}

\begin{proof}
Suppose $G$ has right adjoint $U$, $H$ has right adjoint $V$, and $\Phi
$ has right adjoint $\Gamma $.  The natural transformation $\tau $
induces a dual natural transformation $D\tau \mathcolon \Gamma
V\xrightarrow{}U\Gamma $.  Define $\funcBF{\Gamma } \mathcolon
\genBF{\cat{D}}{H}\xrightarrow{} \BF $ by $(\funcBF{\Gamma }
Y)_{n}=\Gamma Y_{n}$, with structure maps adjoint to the composite
\[
\Gamma Y_{n}\xrightarrow{\Gamma \widetilde{\sigma }}\Gamma VY_{n+1}
\xrightarrow{D\tau } U\Gamma Y_{n+1}
\]
where $\widetilde{\sigma }$ is adjoint to the structure map of $Y$.  The
functor $\funcBF{\Gamma }$ is analogous to restriction in the theory of
group representations, and we must now define the analog to induction
$\funcBF{\Phi }$.  To do so, first note that $\tau $ defines natural
transformations $\tau ^{q}\mathcolon \Phi G^{q}\xrightarrow{}H^{q}\Phi $
for all $q$, by iteration.  Define $(\funcBF{\Phi }X)_{n}$ to be the
coequalizer of the two maps
\[
\coprod _{p+q+r=n} H^{p}\Phi G^{q}X_{r} \rightrightarrows \coprod _{p+q=n}
H^{p}\Phi X_{q} 
\]
where the top map is induced by $H^{p}\Phi
G^{q}X_{r}\xrightarrow{}H^{p}\Phi X_{q+r}$ and the bottom map is induced
by $H^{p}\Phi G^{q}X_{r}\xrightarrow{H^{p}\tau ^{q}}H^{p+q}\Phi X_{r}$.
To define the structure map of $\funcBF{\Phi }X$, note that the
coequalizer diagram for $H(\funcBF{\Phi }X)_{n}$ is just the subdiagram
of the coequalizer diagram for $(\funcBF{\Phi }X)_{n+1}$ consisting of
all terms with a positive power of $H$.  The inclusion of diagrams
induces the desired structure map $H(\funcBF{\Phi
}X)_{n}\xrightarrow{}(\funcBF{\Phi }X)_{n+1}$.  We leave to the reader
the exercise in adjointness required to prove that $\funcBF{\Phi }$ is left
adjoint to $\funcBF{\Gamma }$.

The functor $\funcBF{\Gamma }$ clearly preserves level fibrations and
level trivial fibrations, so $\funcBF{\Phi }$ is a left Quillen functor
with respect to the projective model structures.  Also, since $\Ev
_{n}\funcBF{\Gamma }=\Gamma \Ev _{n}$, we have $F_{n}\Phi =\funcBF{\Phi
}F_{n}$.  To show that $\funcBF{\Phi }$ is a left Quillen functor with
respect to the stable model structures, we must show that $\funcBF{\Phi
}s_{n}^{A}$ is a stable equivalence for all cofibrant $A$, by
Theorem~\ref{thm-localization}.  Using the fact that $\funcBF{\Phi
}F_{n}=F_{n}\Phi $ and the fact that $\tau _{A}$ is a weak equivalence,
we reduce to showing that $s_{n}\Phi A$ is a stable equivalence in
$\genBF{\cat{D}}{H}$.  This follows from Theorem~\ref{thm-stable-equiv},
so $\funcBF{\Phi }$ is a left Quillen functor with respect to the stable
model structures.  the map $\Phi s_{n}^{A}$ is a stable equivalence.
Thus $\Phi \mathcolon \BF \xrightarrow{}\genBF{\cat{D}}{H}$ is a left
Quillen functor with respect to the stable model structures.

We define $\funcBF{\tau }$ by defining its adjoint $D\funcBF{\tau
}\mathcolon \funcBF{\Gamma }V\xrightarrow{}U \funcBF{\Gamma }$.  Indeed,
$D\tau $ is just the prolongation of the adjoint $D\tau $ of $\tau $.
Since $\tau $ is a weak equivalence on all cofibrant objects of
$\cat{C}$, $D\tau $ is a weak equivalence on all fibrant objects of
$\cat{D}$.  To see this, note that $\tau $ induces a natural isomorphism
in the homotopy category.  Adjointness implies that $D\tau $ also
induces a natural isomorphism in the homotopy category, and it follows
that $D\tau $ is a weak equivalence on all fibrant objects of $\cat{D}$.
Thus $D\funcBF{\tau }$ will be a level equivalence on all level fibrant
objects of $\genBF{\cat{D}}{H}$, so $\funcBF{\tau }$ is a level
equivalence on all cofibrant objects of $\BF $.  We leave it to the
reader to check compatibility of $(\funcBF{\Phi },\funcBF{\tau })$ with
compositions and identities.
\end{proof}

Proposition~\ref{prop-functor} and Theorem~\ref{thm-already-equiv} give
us a weak universal property of $\BF $.  

\begin{corollary}\label{cor-weak-universality}
Suppose $(\Phi ,\tau )\mathcolon (\cat{C},G) \xrightarrow{}(\cat{D},H)$
is a map of pairs such that $H$ is a Quillen equivalence.  Then there is
a functor $\widetilde{\Phi }\mathcolon \ho \BF \xrightarrow{}\cat{D}$
such that $\widetilde{\Phi }\circ LF_{0}=L\Phi \mathcolon \ho
\cat{C}\xrightarrow{}\ho \cat{D}$.
\end{corollary}

This corollary is trying to say that $(\BF ,G)$ is homotopy initial among
maps of pairs $(\cat{C},G)\xrightarrow{}(\cat{D},H)$ where $H$ is a
Quillen equivalence.  Though the statement of the corollary is the best
statement of this concept we have been able to find, we suspect there
is a better one.  

\begin{proof}
By Proposition~\ref{prop-functor} there is a map of pairs 
\[
(\funcBF{\Phi },\funcBF{\tau })\mathcolon (\BF,G)
\xrightarrow{}(\genBF{\cat{D}}{H},H)
\]
induced by $(\Phi ,\tau )$.  By Theorem~\ref{thm-already-equiv},
$F_{0}\mathcolon \cat{D}\xrightarrow{}\genBF{\cat{D}}{H}$ is a Quillen
equivalence.  Define $\widetilde{\Phi }$ to be the composite $R\Ev _{0}\circ
L\funcBF{\Phi }$.
\end{proof}

We have now shown that the correspondence $(\cat{C},G)\mapsto (\BF ,G)$
is functorial.  We would like to know that it is homotopy invariant.  In
particular, we would like to know that $(\funcBF{\Phi }, \funcBF{\tau
})$ is a Quillen equivalence of pairs when $\Phi $ is a Quillen
equivalence.  Our proof of this seems to require some hypotheses.  

\begin{theorem}\label{thm-equiv}
Suppose $(\Phi ,\tau )\mathcolon (\cat{C},G) \xrightarrow{}(\cat{D},H)$
is a map of pairs such that $\Phi $ is a Quillen equivalence.  Suppose
as well that either the domains of the generating cofibrations for
$\cat{C}$ can be taken to be cofibrant, or that $\tau_{X}$ is a weak
equivalence for all $X$.  Then, in the induced map of pairs
$(\funcBF{\Phi }, \funcBF{\tau })\mathcolon (\BF ,G) \xrightarrow{}
(\genBF{\cat{D}}{H},H)$, the Quillen functor $\funcBF{\Phi }$ is a
Quillen equivalence. 
\end{theorem}

\begin{proof}
We will first show that $\funcBF{\Phi }$ is a Quillen equivalence on the
projective model structures.  Use the same notation as in the proof of
Proposition~\ref{prop-functor}, so that $\Gamma $ denotes the right
adjoint of $\Phi $.  Then, since $\Phi $ is a Quillen equivalence,
$\Gamma $ reflects weak equivalences between fibrant objects,
by~\cite[Corollary~1.3.16]{hovey-model}.  It follows that
$\funcBF{\Gamma }$ reflects level equivalences between level fibrant
objects.  Hence to show that $\funcBF{\Phi }$ is a Quillen equivalence
of the projective model structures, it suffices to show that
$X\xrightarrow{}\funcBF{\Gamma }R\funcBF{\Phi }X$ is a level equivalence
for all cofibrant $X$, where $R$ is a fibrant replacement functor in the
projective model structure on $\genBF{\cat{D}}{H}$.  Thus, we need only
show that $X_{n}\xrightarrow{}\Gamma R(\funcBF{\Phi }X)_{n}$ is a weak
equivalence for all $n$ and all cofibrant $X$, where now $R$ is a
fibrant replacement functor in $\cat{D}$.  Since $X_{n}$ is cofibrant
and $\Phi $ is a Quillen equivalence, it suffices to show that $\Phi
X_{n}\xrightarrow{}(\funcBF{\Phi }X)_{n}$ is a weak equivalence for all
$n$ and all cofibrant $X$. In fact, we can assume that $X$ is an
$I_{G}$-cell complex.

Write $X$ as the colimit of a $\lambda $-sequence 
\[
0=X^{0}\xrightarrow{}X^{1} \xrightarrow{}X^{2}\xrightarrow{} \dots
\xrightarrow{} X^{\beta } \xrightarrow{} \dots \xrightarrow{} X^{\lambda
}=X 
\]
where each map $X^{\beta }\xrightarrow{}X^{\beta +1}$ is a pushout of a
map of $I_{G}$.  We will prove that, for all $\beta \leq \lambda $,
$\Phi X^{\beta }_{n}\xrightarrow{}(\funcBF{\Phi }X^{\beta })_{n}$ is a
weak equivalence for all $n$, by transfinite induction on $\beta $.
Getting started is easy.  The limit ordinal part of the induction
follows from~\cite[Proposition 17.9.12]{hirschhorn}, since each of the
maps $\Phi X^{\beta }_{n}\xrightarrow{}\Phi X^{\beta +1}_{n}$ and each
of the maps $(\funcBF{\Phi }X^{\beta })_{n}\xrightarrow{}(\funcBF{\Phi
}X^{\beta +1})_{n}$ is a cofibration of cofibrant objects.

For the succesor ordinal part of the induction, suppose $X^{\beta
}\xrightarrow{}X^{\beta +1}$ is a pushout of the map
$F_{m}C\xrightarrow{F_{m}f}F_{m}D$ of $I_{G}$.  Then we have a pushout
diagram 
\[
\begin{CD}
\Phi (F_{m}C)_{n} @>>> \Phi (F_{m}D)_{n} \\
@VVV @VVV \\
\Phi X^{\beta }_{n} @>>> \Phi X^{\beta +1}_{n}
\end{CD}
\]
and another pushout diagram 
\[
\begin{CD}
(\funcBF{\Phi }F_{m}C)_{n} @>>> (\funcBF{\Phi }F_{m}D)_{n} \\
@VVV @VVV \\
(\funcBF{\Phi }X^{\beta })_{n} @>>> (\funcBF{\Phi }X^{\beta +1})_{n}
\end{CD}
\]
in $\cat{D}$.  Note that $\Phi (F_{m}C)_{n}=\Phi G^{n-m}C$, where we
interpret $G^{n-m}C$ to be the initial object if $n<m$.  Similarly,
$(\funcBF{\Phi }F_{m}C)_{n}=(F_{m}\Phi C)_{n}=H^{n-m}\Phi C$.  Thus the
natural transformation $\tau $ induces a map from the first of these
pushout squares to the second.  If $C$ (and hence also $D$) is
cofibrant, then this map of pushout squares is a weak equivalence at
both the upper left and upper right corners.  Or, if $\tau $ is a weak
equivalence for all $X$, then again this map is a weak equivalence at
both the upper left and upper right corners.  It is also a weak
equivalence at the lower left corner, by the induction hypothesis.
Since the top horizontal map is a cofibration in $\cat{D}$, Dan Kan's
cube lemma~\cite[Lemma 5.2.6]{hovey-model} implies that the map is a
weak equivalence on the lower right corner.  This completes the
induction.

We have now proved that $\funcBF{\Phi }$ is a Quillen equivalence with
respect to the projective model structures.  In view of
Proposition~\ref{prop-local-quillen}, to show that $\funcBF{\Phi }$ is a
Quillen equivalence with respect to the stable model structures, we need
to show that if $Y$ is level fibrant in $\genBF{\cat{D}}{H}$ and
$\funcBF{\Gamma }Y$ is a $U$-spectrum, then $Y$ is a $V$-spectrum.  To
see this, note that, since $\funcBF{\Gamma }Y$ is a $U$-spectrum, the
natural map $\Gamma Y_{n}\xrightarrow{}U\Gamma Y_{n+1}$ is a weak
equivalence for all $n$.  There is a natural transformation $D\tau
\mathcolon \Gamma V\xrightarrow{}U\Gamma $ dual to $\tau $.
Furthermore, $(D\tau)_{X}$ is a weak equivalence for all fibrant $X$, as
we have seen in the proof of Proposition~\ref{prop-functor}.  Thus, the
natural map $\Gamma Y_{n}\xrightarrow{}\Gamma VY_{n+1}$ is a weak
equivalence for all $n$.  Since $\Gamma $ reflects weak equivalences
between fibrant objects, it follows that $Y$ is a $V$-spectrum, as
required.
\end{proof}

As an example of Theorem~\ref{thm-equiv}, suppose we take a pointed
simplicial set $K$ weakly equivalent to $S^{1}$.  Then there is a weak
equivalence $K\xrightarrow{}RS^{1}$, where $R$ is the fibrant
replacement functor.  This induces a natural transformation of left
Quillen functors $\tau \mathcolon K\smash -\xrightarrow{}RS^{1}\smash
-$.  In Theorem~\ref{thm-equiv}, take $\cat{D}=\cat{C}$ equal to the
model category of pointed simplicial sets, take $\Phi $ to be the identity,
and take $\tau $ to be this natural transformation.  Then we get a
Quillen equivalence between the stable model categories of spectra
obtained by inverting $K$ and inverting $RS^{1}$.  Therefore, the choice
of simplicial circle does not matter, up to Quillen equivalence, for
Bousfield-Friedlander spectra.  

We now investigate to what extent the correspondence that takes
$(\cat{C},G)$ to the stable model structure on $\BF $ preserves monoidal
structure.  We begin by assuming that $\cat{C}$ is a $\cat{D}$-model
category, for some symmetric monoidal model category $\cat{D}$.  This
means that $\cat{D}$ is a symmetric monoidal category with a compatible
model structure.  Since we will need to work with this compatibility, we
remind the reader of the precise definition (see
also~\cite[Chapter~4]{hovey-model}).  

\begin{definition}\label{defn-boxprod}
Suppose $\cat{D}$ is a monoidal category.  Given maps $f\mathcolon
A\xrightarrow{}B$ and $g\mathcolon C\xrightarrow{}D$, we define the
\emph{pushout product} $f\boxprod g$ of $f$ and $g$ to be the map
\newline $f\boxprod g\mathcolon (A\otimes D)\amalg _{A\otimes C}(B\otimes 
C)\xrightarrow{}B\otimes D$ induced by the commutative square 
\[
\begin{CD}
A\otimes C @>f\otimes 1>> B\otimes C \\
@V1\otimes gVV @VV1\otimes gV \\
A\otimes D @>>f\otimes 1> B\otimes D
\end{CD}
\]
\end{definition}

The compatibility condition we require is then that, if $f$ and $g$ are
cofibrations, then so is $f\boxprod g$, and furthermore, if one of $f$
and $g$ is a trivial cofibration, so is $f\boxprod g$.  We must also
require that, if $S$ is the unit of $\otimes $ and $QS\xrightarrow{}S$
is a cofibrant approximation, then $QS\otimes X\xrightarrow{}X$ is still
a weak equivalence.  

Then, by saying that $\cat{C}$ is a $\cat{D}$-model category, we mean
that $\cat{C}$ is tensored, cotensored, and enriched over $\cat{D}$,
compatibly with the model structure.  This compatibility is precisely
analogous to the compatibility above.  

We then have the following theorem.  

\begin{theorem}\label{thm-spectra-are-simplicial}
Let $\cat{D}$ be a cofibrantly generated monoidal model category, and
suppose the domains of the generating cofibrations are cofibrant.
Suppose $\cat{C}$ is a left proper cellular $\cat{D}$-model category,
and that $G$ is a left $\cat{D}$-Quillen endofunctor of $\cat{D}$.  This
means that $G(X\otimes K)\cong GX\otimes K$, coherently, for $X\in
\cat{C}$ and $K\in \cat{D}$.  Then $\BF $, with the stable model
structure, is again a $\cat{D}$-model category, and the extension of $G$
is a $\cat{C}$-Quillen self-equivalence of $\BF $.
\end{theorem}

Of course, the Quillen functors $F_{n}\mathcolon
\cat{C}\xrightarrow{}\BF$ will be $\cat{D}$-Quillen
functors as well.

\begin{proof}
We define the action of $\cat{D}$ on $\BF $ levelwise.  That is, given
$X\in \BF$ and $K\in \cat{D}$, we define $(X\otimes K)_{n}=X_{n}\otimes
K$.  The structure map is given by
\[
G(X_{n}\smash K)\cong GX_{n}\smash K \xrightarrow{}X_{n+1}\smash K . 
\]
One can easily verify that this makes $\BF $ tensored over $\cat{D}$.
Similarly, define $(X^{K})_{n}=X_{n}^{K}$, with structure maps
$G(X_{n}^{K})\xrightarrow{}X_{n+1}^{K}$ adjoint to the composite
\[
G(X_{n}^{K})\smash K \cong G(X_{n}^{K}\smash K)
\xrightarrow{G(\text{ev})} GX_{n} \xrightarrow{}X_{n+1} 
\]
where $\text{ev}\mathcolon X_{n}^{K}\smash K$ is the evaluation map,
adjoint to the identity of $X_{n}^{K}$.  This makes $\BF $ cotensored
over $\cat{D}$.  Finally, given $X$ and $Y$ in $\BF $, 
define $\Map (X,Y)\in \cat{D}$ to be the equalizer of the
two maps 
\[
\alpha ,\beta \mathcolon \prod _{n}\Map (X_{n},Y_{n})\xrightarrow{}\prod
_{n}\Map (X_{n},UY_{n+1})
\]
where $\alpha $ is the product of the maps $\Map
(X_{n},Y_{n})\xrightarrow{}\Map (X_{n},UY_{n+1})$ induced by the adjoint
of the structure map of $Y$, and $\beta $ is the product of the maps 
\[
\Map (X_{n+1},Y_{n+1}) \xrightarrow{}\Map (UX_{n+1},UY_{n+1})
\xrightarrow{} \Map (X_{n},UY_{n+1}).
\]
Here the first map exists since $G$ preserves the $\cat{D}$ action, and
the second map is induced by the structure map of $X$.  This functor
makes $\BF $ enriched over $\cat{D}$.  

We must now check that these structures are compatible with the model
structure.  We begin with the projective model structure on $\BF $.  One
can easily check that $F_{n}f\boxprod g = F_{n}(f\boxprod g)$.  Thus, if
$f$ is one of the generating cofibrations of the projective model
structure on $\BF $, and $g$ is a cofibration in $\cat{D}$, then
$f\boxprod g$ is a cofibration in $\BF $.  It follows that
$f\boxprod g$ is a cofibration for $f$ an arbitrary cofibration of
$\BF $ (see~\cite[Lemma 2.3]{schwede-shipley-monoids}
and~\cite[Corollary 5.3.5]{hovey-shipley-smith}).  A similar argument
shows that $f\boxprod g$ is a projective trivial cofibration in $\BF $ if
either $f$ is a projective cofibration in $\BF$ or $g$ is a trivial
cofibration in $\cat{D}$.  Finally, if $QS\xrightarrow{}S$ is a
cofibrant approximation to the unit $S$ in $\cat{D}$, and $X$ is
cofibrant in $\BF $, then each 
$X_{n}$ is cofibrant in $\cat{C}$, so the map $X\otimes 
QS\xrightarrow{}X$ is a level equivalence as required.  Thus
$\BF $ with its projective model structure is a $\cat{D}$-model category.

To show that $\BF$ with its stable model structure is also a
$\cat{D}$-model category, we need to show that, if $f$ is a stable
trivial cofibration and $g$ is a cofibration in $\cat{D}$, then
$f\boxprod g$ is a stable equivalence.  It suffices to check this for
$g\mathcolon K\xrightarrow{}L$ one of the generating trivial
cofibrations of $\cat{D}$.  In this case, by hypothesis, $K$ and $L$ are
cofibrant in $\cat{D}$.  Thus the functor $-\otimes K$ is a Quillen
functor with respect to the projective model structure on $\BF $, and
similarly for $L$.  Furthermore, if $s_{n}^{QC}\mathcolon
F_{n+1}GQC\xrightarrow{}F_{n}QC$ is an element of the set $\cat{S}$, then
$s_{n}^{QC}\otimes K\cong s_{n}^{QC\otimes K}$, since $G$ preserves the
$\cat{D}$-action.  In view of Theorem~\ref{thm-stable-equiv}, the map
$s_{n}^{QC\otimes K}$ is a stable equivalence.
Theorem~\ref{thm-localization} then implies that $-\otimes K$ is a
Quillen functor with respect to the stable model structure on $\BF $,
and similarly for $-\otimes L$.  Thus, if $f$ is a stable trivial
cofibration, so are $f\otimes K$ and $f\otimes L$.  It follows from the
two out of three property that $f\boxprod g$ is a stable equivalence, as
required.
\end{proof}

\begin{remark}\label{rem-subtle-action}
Suppose that the functor $G$ is actually given by $GX=X\otimes K$ for
some cofibrant object $K$ of $\cat{D}$.  We then have \textbf{two
different ways} of tensoring with $K$ on $\BF $.  The first way is the
extension of $G$ to a Quillen equivalence of $\BF $.  Recall from
Remark~\ref{rem-subtle} that this functor, which we denote by $X\mapsto
X\stensor K$, does not use the twist map.  On the other hand, we also
have the functor $X\mapsto X\otimes K$ that is part of the
$\cat{D}$-action on $\BF $ constructed in
Theorem~\ref{thm-spectra-are-simplicial}.  This functor \emph{does} use
the twist map as part of its structure map; indeed, in order to
construct the isomorphism $G(X\otimes K)\cong GX\otimes K$ we need to
permute the two different copies of $K$.  Therefore, we \textbf{do not
know} that $X\mapsto X\otimes K$ is a Quillen equivalence, even though
$X\mapsto X\stensor K$ is.  We will have to deal with this point more
thoroughly in Section~\ref{sec-comparison}, when we compare $\BF $ with
symmetric spectra.  
\end{remark}

Theorem~\ref{thm-spectra-are-simplicial} gives us a functorial
stabilization.  We first simplify the notation.  Suppose $K$ is a
cofibrant object of a symmetric monoidal model category $\cat{D}$.  Then
$G= -\otimes K$ is a left Quillen functor on any $\cat{D}$-model
category $\cat{C}$.  In this case, we denote $\BF $ by
$\genBF{\cat{C}}{K}$.  

\begin{corollary}\label{cor-functorial-stabilization}
Suppose $K$ is a cofibrant object of a cofibrantly generated symmetric
monoidal model category $\cat{D}$ where the domains of the generating
cofibrations can be taken to be cofibrant.  Then the correspondence
$\cat{C}\mapsto \genBF{\cat{C}}{K} $ defines an endofunctor of the
category of left proper cellular $\cat{D}$-model categories.
\end{corollary}

Note that the ``category'' of left proper cellular $\cat{D}$-model
categories is not really a category, because the $\Hom $-sets need not
be sets.  It is really a $2$-category, and the correspondence
$\cat{C}\mapsto \genBF{\cat{C}}{K}$ is actually a $2$-functor.
See~\cite{hovey-model} for a description of this point of view on model
categories.

\begin{proof}
Given a left proper cellular $\cat{D}$-model category
$\cat{C}$, we have seen in Theorem~\ref{thm-spectra-are-simplicial} that
$\genBF{\cat{C}}{K} $ is a $\cat{D}$-model category.  Just as in
Lemma~\ref{lem-extending-functors}, a $\cat{D}$-Quillen functor
$H\mathcolon \cat{C}\xrightarrow{}\cat{C}'$ induces a functor
$H\mathcolon \genBF{\cat{C}}{K}\xrightarrow{}\genBF{\cat{C}'}{K}$, as
does its right adjoint $V$.  Since $H$ is defined levelwise, it
preserves the action of $\cat{D}$.  It is easy to check that $V$
preserves level fibrations and level trivial fibrations, so that $H$ is
a $\cat{C}$-Quillen functor with respect to the projective model
structures.  Furthermore, we have
$Hs_{n}^{QC}=s_{n}^{HQC}$, so Theorem~\ref{thm-stable-equiv} and
Theorem~\ref{thm-localization} imply that $H$ is a $\cat{C}$-Quillen
functor with respect to the stable model structures as well.  
\end{proof}

We now point out that, if $\cat{D}$ is a symmetric monoidal model
category, and $G=-\otimes K$ for some cofibrant object $K$, the category
$\genBF{\cat{D}}{G}$ is almost never itself monoidal, though, as we have
seen, it has an action of $\cat{D}$.  To see this, consider the category
$\cat{D}^{\mathbb{N}}$ of sequences from $\cat{D}$.  An object of
$\cat{D}^{\mathbb{N}}$ is a sequence $X_{n}$ of objects of $\cat{D}$,
and a map $f\mathcolon X\xrightarrow{}Y$ is a sequence of maps
$f_{n}\mathcolon X_{n}\xrightarrow{}Y_{n}$.  Then $\cat{D}^{\mathbb{N}}$
is a symmetric monoidal category, where we define $(X\smash
Y)_{n}=\coprod _{p+q=n}X_{p}\smash Y_{q}$.  The functor $G$ defines a
monoid $T=(S^{0},GS^{0},G^{2}S^{0},\dots ,G^{n}S^{0},\dots )$ in this
category, using the fact that $G$ preserves the $\cat{D}$-action. 

\begin{lemma}\label{lem-modules-over-G}
Suppose $\cat{D}$ is a symmetric monoidal model category and $G$ is a left
$\cat{D}$-Quillen functor.  Then $\genBF{\cat{D}}{G}$ is the category of left
modules over the monoid $T=(S^{0},GS^{0},\dots ,G^{n}S^{0},\dots )$.  
\end{lemma}

We leave the proof of this lemma to the reader.  The important point is
that the monoid $T$ is almost never commutative, and therefore
$\genBF{\cat{D}}{G}$ can not be a symmetric monoidal category with unit
$T$.  Indeed, let $K=GS\in \cat{D}$, so that $G=-\otimes K$.  Then $T$
is commutative if and only if the commutativity isomorphism on $K\smash
K$ is the identity.  This happens only very rarely.

\section{Symmetric spectra}\label{sec-symm-spec}

We have just seen that the stabilization functor $\BF $ is not good
enough in case $\cat{D}$ is a symmetric monoidal model category and $G$
is a $\cat{D}$-Quillen functor, because $\BF $ is not usually itself a
symmetric monoidal model category.  In this section, we begin the
construction of a better stabilization functor $\spec $ for
this case.  We will concentrate on the category theory in this section,
leaving the model structures for the next section.  The terms used for
the algebra of symmetric monoidal categories and modules over them are
all defined in~\cite[Section~4.1]{hovey-model}.  

Through most of this section, then, $\cat{D}$ will be a bicomplete
closed symmetric monoidal category with unit $S$, and $K$ will be an
object of $\cat{D}$, The category $\cat{C}$ will be a bicomplete
category enriched, tensored, and cotensored over $\cat{D}$.  Note that
any $\cat{D}$-functor $G$ on $\cat{D}$ is of the form $G(X)=X\smash K$
for $K=GS$, so we will only consider such functors.  Because of this, we
will drop the letter $G$ from our notations and replace it with $K$.

This section is based on the symmetric spectra and sequences
of~\cite{hovey-shipley-smith}.  The main idea of symmetric spectra is
that the commutativity isomorphism of $\cat{D}$ makes $K^{\smash n}$
into a $\Sigma _{n}$-object of $\cat{D}$, where $\Sigma _{n}$ is the
symmetric group on $n$ letters.  We must keep track of this action if we
expect to get a symmetric monoidal category of $K$-spectra.

The following definition is~\cite[Definition
2.1.1]{hovey-shipley-smith}.

\begin{definition}\label{defn-symmetric-sequence}
Let $\Sigma =\coprod _{n\geq 0}\Sigma _{n}$ be the category whose
objects are the sets $\overline{n}=\{1,2,\dots ,n \}$ for $n\geq 0$,
where $\overline{0}=\emptyset $.  The morphisms of $\Sigma $ are the
isomorphisms of $\overline{n}$.  Given a category $\cat{C}$, a
\emph{symmetric sequence} in $\cat{C}$ is a functor $\Sigma
\xrightarrow{}\cat{C}$.  The category of symmetric sequences is the
functor category $\cat{C}^{\Sigma }$.
\end{definition}

A symmetric sequence in $\cat{C}$ is a sequence $X_{0},X_{1},\dots
,X_{n},\dots $ of objects of $\cat{C}$ with an action of $\Sigma _{n}$
on $X_{n}$.  It is sometimes more useful to consider a symmetric
sequence as a functor from the category of finite sets and isomorphisms
to $\cat{C}$; since the category $\Sigma $ is a skeleton of the category
of finite sets and isomorphisms, there is no difficulty in doing so.  

As a functor category, the category of symmetric sequences in $\cat{C}$
is bicomplete if $\cat{C}$ is so; limits and colimits are taken
objectwise.  Furthermore, if $\cat{D}$ is a closed symmetric monoidal
category, so is $\cat{D}^{\Sigma }$, as explained
in~\cite[Section~2.1]{hovey-shipley-smith}.  Recall that the monoidal
structure is given by 
\[
(X\otimes Y)(C) = \coprod _{A\cup B=C,A\cap B=\emptyset } X(A)\otimes Y(B)
\]
where we think of $X$, $Y$, and $X\otimes Y$ as functors from finite sets
to $\cat{D}$.  Equivalently, though less canonically, we have 
\[
(X\otimes Y)_{n} = \coprod _{p+q=n} \Sigma _{n}\times _{\Sigma
_{p}\times \Sigma _{q}} (X_{p}\otimes Y_{q}).
\]
The unit of the monoidal structure is the symmetric
sequence $(S,0,\dots ,0,\dots )$, where $0$ is the initial object of
$\cat{C}$.  To define the closed structure, we first define $\Hom
_{\Sigma _{n}}(X,Y)$ for $X,Y\in \cat{D}^{\Sigma _{n}}$ in the usual
way, as an equalizer of the two obvious maps $\Hom (X,Y) \xrightarrow{}
\Hom (X\times \Sigma _{n},Y)$.  The closed structure is then given by
\[
\Hom (X,Y)_{k}=\prod _{n}\Hom_{\Sigma _{n}} (X_{n},Y_{n+k}).
\]

If $\cat{C}$ is enriched, tensored, and cotensored over $\cat{D}$, then
$\cat{C}^{\Sigma }$ is enriched, tensored, and cotensored over
$\cat{D}^{\Sigma }$.  Indeed, the same 
definition as above works to define the tensor structure.  The cotensor
structure is defined as follows.  First we define $\Hom _{\Sigma _{n}}(K,X)$
for $X\in \cat{C}^{\Sigma _{n}}$ and $K\in \cat{D}^{\Sigma _{n}}$ as
an appropriate equalizer.  Then, for $X\in \cat{C}^{\Sigma }$ and $K\in
\cat{D}^{\Sigma }$, we define $X^{K}_{k}=\prod _{n}\Hom _{\Sigma
_{n}}(K_{n},X_{n+k})$.  The enrichment $\Map (X,Y)$ is defined
similarly.  In the same way, if $\cat{C}$ is an enriched monoidal
category over $\cat{D}$, then $\cat{C}^{\Sigma }$ is an enriched
monoidal category over $\cat{D}^{\Sigma }$.  

Consider the free commutative monoid $\Sym (K)$ on the object
$(0,K,0,\dots ,0,\dots )$ of $\cat{D}^{\Sigma }$.  One can easily check
that $\Sym (K)$ is the symmetric sequence $(S^{0},K,K\otimes K,\dots
,K^{\otimes n},\dots )$ where $\Sigma _{n}$ acts on $K^{\otimes n}$ by the
commutativity isomorphism, as
in~\cite[Section~4.4]{hovey-shipley-smith}.  

\begin{definition}\label{defn-symm-spectra}
Suppose $\cat{D}$ is a symmetric monoidal model category, $\cat{C}$ is a
$\cat{D}$-model category, and $K$ is an object of $\cat{D}$.  The
category of \emph{symmetric spectra} $\modspec $ is the
category of modules in $\cat{C}^{\Sigma }$ over the commutative monoid
$\Sym (K)$ in $\cat{D}^{\Sigma }$.  That is, a symmetric spectrum is a
sequence of objects $X_{n}\in \cat{C}^{\Sigma _{n}}$ and $\Sigma
_{n}$-equivariant maps $X_{n}\otimes K \xrightarrow{}X_{n+1}$, such that
the composite 
\[
X_{n} \otimes K^{\otimes p} \xrightarrow{} X_{n+1} \otimes K^{\otimes p-1}
\xrightarrow{} \dots \xrightarrow{} X_{n+p}
\]
is $\Sigma _{n}\times \Sigma _{p}$-equivariant for all $n,p\geq 0$.  A
map of symmetric spectra is a collection of $\Sigma _{n}$-equivariant
maps $X_{n}\xrightarrow{}Y_{n}$ compatible with the structure maps of
$X$ and $Y$.
\end{definition}

Because $\Sym (K)$ is a commutative monoid, the category
$\spec $ is a bicomplete closed symmetric monoidal
category, with $\Sym (K)$ itself as the unit (see Lemma~2.2.2 and
Lemma~2.2.8 of~\cite{hovey-shipley-smith}).  We denote the monoidal
structure by $X\smash Y = X\otimes _{\Sym (K)}Y$, and the closed
structure by $\Hom _{\Sym (K)}(X,Y)$.  Similarly,
$\modspec $ is bicomplete, enriched, tensored, and cotensored over
$\spec $ with the tensor structure denoted $X\smash Y$ again, and, if 
$\cat{C}$ is a $\cat{D}$-monoidal model category, then $\modspec $ will
be a monoidal category enriched over $\spec $. 

Of course, if we take $\cat{C}=\sset _{*}$ and $K=S^{1}$, we recover the
definition of symmetric spectra given in~\cite{hovey-shipley-smith},
except that we are using right modules instead of left modules.  

\begin{definition}\label{defn-symm-ev}
Given $n\geq 0$, the \emph{evaluation functor} $\Ev _{n}\mathcolon
\modspec \xrightarrow{} \cat{C}$ takes $X$ to $X_{n}$.  the evaluation
functor has a left adjoint $F_{n}\mathcolon
\cat{C}\xrightarrow{}\modspec $, defined by $F_{n}X=\widetilde{F}_{n}X
\otimes \Sym (K)$, where $\widetilde{F}_{n}X$ is the symmetric sequence
$(0,\dots ,0,\Sigma _{n}\times X,0,\dots )$.
\end{definition}

Note that $F_{0}X=(X,X\otimes K,\dots ,X\otimes K^{\otimes n},\dots )$,
and in particular $F_{0}S=\Sym (K)$.  Also, if $X\in \cat{C}$ and $Y\in
\cat{D}$, there is a natural isomorphism $F_{n}X\smash F_{m}Y\cong
F_{n+m}(X\otimes Y)$, just as
in~\cite[Proposition~2.2.6]{hovey-shipley-smith}.  In particular,
$F_{0}\mathcolon \cat{D}\xrightarrow{}\spec $ is a (symmetric) monoidal
functor, and so $\modspec $ is naturally enriched, tensored, and
cotensored over $\cat{D}$.  In fact, this structure is very simple.
Indeed, if $X\in \modspec $ and $L\in \cat{D}$, $X\otimes L=X\otimes
_{\Sym (K)}F_{0}L$ is just the symmetric sequence whose $n$th 
term is $X_{n}\otimes L$.  The structure map is the composite
\[
X_{n}\otimes L\otimes K \xrightarrow{1\otimes \tau } X_{n}\otimes K\otimes L
\xrightarrow{} X_{n+1}\otimes L .
\]
Note the presence of the twist map; this is required even when $L=K$ to
get a symmetric spectrum, unlike the case of ordinary spectra.  
Similarly, $X^{L}=\Hom _{\Sym (K)}(F_{0}L,X)$ is the symmetric sequence
whose $n$th term is $X_{n}^{L}$, with the twist map again appearing as
part of the structure map.

\begin{remark}\label{rem-symm-R}
Just as in the spectrum case, the functor $\Ev _{n}$ has a right adjoint
$R_{n}\mathcolon \cat{C} \xrightarrow{} \modspec $.  Indeed,
$R_{n}X=\Hom (\Sym (K),\widetilde{R}_{n}L)$, where $\widetilde{R}_{n}L$
is the symmetric sequence concentrated in degree $n$ whose $n$th term is
$X^{\Sigma _{n}}$.  
\end{remark}

\section{Model structures on symmetric
spectra}\label{sec-symmetric-model}

Throughout this section, $\cat{D}$ will denote a left proper cellular
symmetric monoidal model category, $\cat{C}$ will denote a left proper
cellular $\cat{D}$-model category, and $K$ will denote a cofibrant
object of $\cat{D}$.  In this section, we discuss the projective and
stable model structures on the category $\modspec $ of symmetric
spectra.  The results in this section are very similar to the
corresponding results for spectra, so we will leave most of the proofs
to the reader.

\begin{definition}
A map $f\in \modspec $ is a \emph{level equivalence} if each
map $f_{n}$ is a weak equivalence in $\cat{C}$.  Similarly, $f$ is a
\emph{level fibration} (resp. \emph{level cofibration}, \emph{level
trivial fibration}, \emph{level trivial cofibration}) if each map
$f_{n}$ is a fibration (resp. cofibration, trivial fibration, trivial
cofibration) in $\cat{C}$.  The map $f$ is a \emph{projective
cofibration} if $f$ has the \llp every level trivial fibration.  
\end{definition}

Then, just as in Definition~\ref{defn-IG}, if we denote the generating
cofibrations of $\cat{C}$ by $I$ and the generating trivial cofibrations
by $J$, we define $I_{K}=\bigcup _{n}F_{n}I$ and $J_{K}=\bigcup
_{n}F_{n}J$.

We have analogues of~\ref{prop-smallness-in-G-spec}--\ref{cor-IG-small}
with the same proofs.  This gives us the projective model structure. 

\begin{theorem}\label{thm-symm-strict}
The projective cofibrations, the level fibrations, and the level
equivalences define a left proper cellular model structure on
$\modspec $.  
\end{theorem}

The set $I_{K}$ is the set of generating cofibrations of the projective
model structure, and $J_{K}$ is the set of generating trivial
cofibrations.  The cellularity of the projective model structure is
proved in the Appendix.

Note that $\Ev _{n}$ takes level (trivial) fibrations to (trivial)
fibrations, so $\Ev _{n}$ is a right Quillen functor and $F_{n}$ is a
left Quillen functor.  

\begin{theorem}\label{thm-symm-proj-mon}
The category $\spec $, with the projective model structure, is a
symmetric monoidal model category.  The category $\modspec $, with its
projective model structure, is a $\spec $-model category.
\end{theorem}

\begin{proof}
We first show that the pushout product $f\boxprod g$ is a (trivial)
cofibration when $f$ is a cofibration in $\modspec $, and $g$ is a
cofibration in $\spec $ (and one of them is a level equivalence).  As
explained in~\cite[Corollary~2.5]{hovey-model}. we may as well assume
that $f$ and $g$ belong to the sets of generating cofibrations or
generating trivial cofibrations.  In either case, we have $f=F_{m}f'$
and $g=F_{n}g'$.  But then $f\boxprod g=F_{m+n}(f'\boxprod g')$.  Since
$F_{m+n}$ is a Quillen functor, the result follows.

Now let $QS$ denote a cofibrant replacement for the unit $S$ in
$\cat{D}$.  Then $F_{0}QS$ is a cofibrant replacement for $F_{0}S=\Sym
(K)$ in $\spec $.  Indeed, $F_{0}QS$ is cofibrant, and $\Ev _{n}F_{0}QS$
is just $QS \otimes K^{\otimes n}$.  Since $K$ is cofibrant and
$\cat{D}$ is a monoidal model category, the map
$F_{0}QS\xrightarrow{}F_{0}S$ is a level equivalence.  Now, if
$X$ is cofibrant in $\modspec $, then each $X_{n}$ is
cofibrant.  Hence the map $X_{n}\otimes QS\xrightarrow{}X_{n}$ is a
weak equivalence for all $n$, and so the map $X\smash
F_{0}QS\xrightarrow{}X$ is a level equivalence, as required.
\end{proof}

We point out here that one can show that the projective model structure
on $\spec $ satisfies the monoid axiom
of~\cite{schwede-shipley-monoids}, assuming that $\cat{D}$ itself does
so.  This means there is a projective model structure on the category of
monoids in $\spec $ and on the category of modules over any monoid.
We do not include the proofs of these statement since we have not been
able to prove the analogous statements for the stable model structure. 

The projective cofibrations of symmetric spectra are more complicated
than they are in the case of ordinary spectra.

\begin{definition}\label{defn-latching}
Define the symmetric spectrum $\overline{\Sym (K)}$ in $\spec $ to be
$0$ in degree $0$ and $K^{\otimes n}$ in degree $n$, for 
$n>0$, with the obvious structure maps.  Define the $n$th \emph{latching
space} $L_{n}X$ of $X\in \modspec $ by $L_{n}X=\Ev
_{n}(X\smash \overline{\Sym (K)})$.
\end{definition}

The obvious map $i\mathcolon \overline{\Sym (K)}\xrightarrow{}\Sym (K)$
induces a $\Sigma _{n}$-equivariant natural transformation $L_{n}X
\xrightarrow{}X$.

Note that the latching space is a $\Sigma _{n}$-object of $\cat{C}$.
There is a model structure on $\Sigma _{n}$-objects of $\cat{C}$ where
the fibrations and weak equivalences are the underlying ones.  This
model structure is cofibrantly generated: if $I$ is the set of
generating cofibrations of $\cat{C}$, then the set of generating
cofibrations of $\cat{C}^{\Sigma _{n}}$ is the set $\Sigma _{n}\times
I$.  Here, for an object $A$, $\Sigma _{n}\times A$ is the coproduct of
$n!$ copies of $A$, given the obvious $\Sigma _{n}$-structure.  

\begin{proposition}\label{prop-symm-cofib}
A map $f\mathcolon X\xrightarrow{}Y$ in $\modspec $ is a projective \ulp
trivial\urp{} cofibration if and only if the induced map $\Ev
_{n}(f\boxprod i)\mathcolon X_{n}\amalg _{L_{n}X}L_{n}Y
\xrightarrow{}Y_{n}$ is a \ulp trivial\urp{} cofibration in
$\cat{C}^{\Sigma _{n}}$ for all $n$.
\end{proposition}

\begin{proof}
We only prove the cofibration case, as the trivial cofibration case is
analogous.  If each map $X_{n}\amalg _{L_{n}X}L_{n}Y
\xrightarrow{}Y_{n}$ is a cofibration, then we can show that $f$ is a
projective cofibration by showing $f$ has the \llp level trivial
fibrations.  Indeed, we construct a lift by induction, just as in the
proof of Proposition~\ref{prop-char-of-proj-cofib}.  To prove the
converse, it suffices to show that $\Ev _{n}(f\boxprod i)$ is a $\Sigma
_{n}$-cofibration for $f\in I_{K}$, since $\Ev _{n}$ is itself a left
Quillen functor.  Then we can write $f=F_{m}g$, and we find that $\Ev
_{n}(f\boxprod i)$ is an isomorphism when $n\neq m$, and is the map
$\Sigma _{m}\times g$ when $n=m$.  This is a $\Sigma _{m}$-cofibration,
as required.
\end{proof}

We must now localize the projective model structure to obtain the stable
model structure.  

\begin{definition}\label{defn-symm-omega-spec}
A symmetric spectrum $X\in \modspec $ is an \emph{$\Omega $-spectrum} if
$X$ is level fibrant and the adjoint $X_{n}\xrightarrow{}X_{n+1}^{K}$ of
the structure map is a weak equivalence for all $n$.
\end{definition}

Just as with Bousfield-Friedlander spectra, we would like the $\Omega
$-spectra to be the fibrant objects in the stable model structure.  We
invert the same maps we did in that case. 

\begin{definition}\label{defn-symm-S}
Define the set of maps $\cat{S}$ in $\modspec $ to be
$\{F_{n+1}(QC\otimes K) \xrightarrow{s_{n}^{QC}} F_{n}QC \}$ as $C$ runs
through the domains and codomains of the generating cofibrations of
$\cat{C}$, and $n\geq 0$.  The map $s_{n}^{QC}$ is adjoint to the map
$QC\otimes K \xrightarrow{} \Ev _{n+1}F_{n}QC=\Sigma _{n+1}\times
(QC\otimes K)$ corresponding to the identity of $\Sigma _{n+1}$.  Define
the \emph{stable model structure} on $\modspec $ to be the Bousfield
localization with respect to $\cat{S}$ of the projective model structure
on $\modspec $.  The $\cat{S}$-local weak equivalences are called the
\emph{stable equivalences}, and the $\cat{S}$-local fibrations are
called the \emph{stable fibrations}.
\end{definition}

The following theorem is then analogous to Theorem~\ref{thm-stable-equiv},
and has the same proof.  

\begin{theorem}\label{thm-symm-stable-equiv}
The stably fibrant symmetric spectra are the $\Omega $-spectra.
Furthermore, for all cofibrant $A\in \cat{C}$ and for all $n\geq 0$, the
map $F_{n+1}(A\otimes K)\xrightarrow{s_{n}^{A}}F_{n}A$ is a stable
equivalence.  
\end{theorem}

Just as in Corollary~\ref{cor-stable-spectra}, this theorem implies
that, when $\cat{C}=\sset _{*}$ or $\top _{*}$, 
$\spec $ is the same as the stable model category on (simplicial or
topological) symmetric spectra discussed in~\cite{hovey-shipley-smith}.

The analog of Corollary~\ref{cor-G-Quillen} also holds, with the same
proof, so that tensoring with $K$ is a Quillen endofunctor of
$\modspec $.  Of course, we want this functor to be a Quillen
equivalence.  As in Definition~\ref{defn-shift}, we prove this by
introducing the shift functors.

\begin{definition}\label{defn-symm-shift}
Define the \emph{right shift functor} $s\mathcolon \modspec
\xrightarrow{}\modspec $ by $sX=\Hom _{\Sym (K)}(F_{1}S,X)$.  Thus
$(sX)_{n}=X_{n+1}$, where the $\Sigma _{n}$-action on $X_{n+1}$ is
induced by the usual inclusion $\Sigma _{n}\xrightarrow{}\Sigma _{n+1}$.
The structure maps of $sX$ are the same as the structure maps of $X$.
Define the \emph{left shift functor} $t\mathcolon \modspec
\xrightarrow{}\modspec $ by $tX=X\otimes _{\Sym (K)}F_{1}S$, so that
$t$ is left adjoint of $s$.  We have $(tX)_{n}=\Sigma _{n}\times
_{\Sigma _{n-1}}X_{n-1}$, with the induced structure maps.
\end{definition}

Note that adjointness gives natural isomorphisms 
\[
(sX)^{K}\cong \Hom _{\Sym (K)}(F_{1}K,X)\cong s(X^{K}).
\]
There is a map $F_{1}K\xrightarrow{}F_{0}S$ which is the identity in
degree $1$.  By adjointness, this map induces a map
\[
X=\Hom _{\Sym (K)}(F_{0}S,X)\xrightarrow{}\Hom _{\Sym
(K)}(F_{1}K,X)=(sX)^{K} .
\]
$X$ is an $\Omega $-spectrum if and only if this map is a level
equivalence and $X$ is level fibrant.  Therefore, the same method used
to prove Theorem~\ref{thm-G-Quillen-equiv} also proves the following
theorem.  

\begin{theorem}\label{thm-symm-Quillen-equiv}
The functors $X\mapsto X\otimes K$ and $t$ are Quillen equivalences on
$\modspec $.  Furthermore, $Rs$ is naturally isomorphic to
$L(-\otimes K)$ and $R((-)^{K})$ is naturally isomorphic to $Lt$.  
\end{theorem}

We have now shown that $\modspec $ is a $K$-stabilization of $\cat{C}$.
However, for this construction to be better than ordinary spectra, we
must show that $\spec $ is a symmetric monoidal model category.

\begin{theorem}\label{thm-symm-monoid}
Suppose that the domains of the generating cofibrations of both
$\cat{C}$ and $\cat{D}$ are cofibrant.  Then the category $\spec $ is a
symmetric monoidal model category, and the category
$\modspec $ is a $\spec $-model category.
\end{theorem}

\begin{proof}
We prove this theorem in the same way as
Theorem~\ref{thm-spectra-are-simplicial}.  Since the cofibrations in the
stable model structure are the same as the cofibrations in the
projective model structure, the only thing to check is that $f\boxprod
g$ is a stable equivalence when $f$ and $g$ are cofibrations and one of
them is a stable equivalence.  We may as well assume that $f\mathcolon
F_{n}A\xrightarrow{}F_{n}B$ is a generating cofibration in $\modspec $
and $g$ is a stable trivial cofibration in $\spec $; the argument for
$f$ a stable trivial cofibration and $g$ a generating cofibration in
$\spec $ is the same.  Then, by hypothesis, $A$ and $B$ are cofibrant in
$\cat{C}$.  We claim that $F_{n}A\smash -$ is a Quillen functor $\spec$
to $\modspec $ with their stable model structures, and similarly for
$F_{n}B\smash -$.  Indeed, in view of Theorem~\ref{thm-localization}, it
suffices to show that $F_{n}A\smash s_{m}^{QC}$ is a stable equivalence
for all $m\geq 0$ and all domains or codomains $C$ of the generating
cofibrations of $\cat{C}$.  But one can easily check that $F_{n}A\smash
s_{m}^{QC}=s_{n+m}^{A\otimes QC}$.  Then
Theorem~\ref{thm-symm-stable-equiv} implies that this map is a stable
equivalence, as required.

Thus, both functors $F_{n}A\smash -$ and $F_{n}B\smash -$ are Quillen
functors in the stable model structures.  A two out of three
argument, as in Theorem~\ref{thm-spectra-are-simplicial}, then shows
that $f\boxprod g$ is a stable equivalence, as required.
\end{proof}

Note that the functor $F_{0}\mathcolon \cat{D}\xrightarrow{}\spec $ is a
symmetric monoidal Quillen functor, so of course $\modspec $ is a
$\cat{D}$-model category as well, under the hypotheses of
Theorem~\ref{thm-symm-monoid}.  In fact, we only need the domains of the
generating cofibrations of $\cat{D}$ to be cofibrant to conclude that
$\modspec $ is a $\cat{D}$-model category, using the argument of
Theorem~\ref{thm-symm-monoid}.  

As we mentioned above, we do not know if the stable model structure
satisfies the monoid axiom.  Given a particular monoid $R$, one could
attempt to localize the projective model structure on $R$-modules to
obtain a stable model structure.  However, for this to work one would
need to know that the projective model structure is cellular, and the
author does not see how to prove this.  This plan will certainly fail
for the category of monoids, since the projective model structure on
monoids will not be left proper in general.  

We also point out that it may be possible to prove some of the results
of Section~\ref{sec-fin-gen} for symmetric spectra.  All of those
results cannot hold, since stable homotopy isomorphisms do not coincide
with stable equivalences even for symmetric spectra of simplicial sets.
Nevertheless, in that case, every stable homotopy isomorphism is a
stable equivalence~\cite[Theorem~3.1.11]{hovey-shipley-smith}, and there
is a replacement for the functor $R^{\infty }$ constructed
in~\cite{shipley-thh}.  We do not know if these results hold for
symmetric spectra over a general well-behaved finitely generated model
category.

\section{Properties of symmetric spectra}\label{sec-symmetric-functor}
 
In this section, we point out that the arguments of
Section~\ref{sec-funct} also apply to symmetric spectra. In particular,
if smashing with $K$ is already a Quillen equivalence on $\cat{C}$, then
$F_{0}\mathcolon \cat{C}\xrightarrow{}\modspec $ is a Quillen
equivalence.  This means that, under mild hypotheses, the homotopy
category of $\cat{C}$ is enriched, tensored, and cotensored over $\ho
\spec $.  We also show symmetric spectra are functorial in an
appropriate sense.  In particular, we show that the Quillen equivalence
class of $\modspec $ is an invariant of the homotopy type of $K$.  

Throughout this section, $\cat{D}$ will denote a left proper cellular
symmetric monoidal model category, $\cat{C}$ will denote a left proper
cellular $\cat{C}$-model category, and $K$ will denote a cofibrant
object of $\cat{D}$.  

The proof of the following theorem is the same as the proof of
Theorem~\ref{thm-already-equiv}.  

\begin{theorem}\label{thm-symm-already-equiv}
Suppose smashing with $K$ is a Quillen equivalence on $\cat{C}$. Then
$F_{0}\mathcolon \cat{C} \xrightarrow{}\modspec $ is a Quillen
equivalence.
\end{theorem}

\begin{corollary}\label{cor-symm-module}
Suppose that the domains of the generating cofibrations of both
$\cat{C}$ and $\cat{D}$ are cofibrant, and suppose that smashing with
$K$ is a Quillen equivalence on $\cat{C}$.  Then $\ho \cat{D}$ is
enriched, tensored, and cotensored over $\ho \spec $.  
\end{corollary}

\begin{proof}
Note that $\ho \modspec $ is certainly enriched, tensored, and
cotensored over $\ho \spec $.  Now use the equivalence of categories
$LF_{0}\mathcolon \ho \cat{C}\xrightarrow{}\ho \modspec $ to transport
this structure back to $\cat{C}$.
\end{proof}

Recall that the homotopy category of any model category is naturally
enriched, tensored, and cotensored over $\ho \sset
$~\cite[Chapter~5]{hovey-model}.  This corollary is the first step to
the assertion that the homotopy category of any stable (with respect to
the suspension) model category is naturally enriched, tensored, and
cotensored over the homotopy category of (simplicial) symmetric
spectra.  See~\cite{schwede-shipley-unique} for further results along
these lines.  

Symmetric spectra are functorial in a natural way.  

\begin{theorem}\label{thm-symm-functorial}
Suppose the domains of the generating cofibrations of $\cat{D}$,
$\cat{C}$, and the left proper cellular $\cat{C}$-model category
$\cat{C}'$ are cofibrant.  Then any $\cat{D}$-Quillen functor
$\Phi \mathcolon \cat{C}\xrightarrow{}\cat{C}'$ extends naturally to a
$\spec $-Quillen functor
\[
\funcspec \mathcolon \modspec \xrightarrow{}\genspec{\cat{C}'}{K}.
\]
Furthermore, if $\Phi $ is a Quillen equivalence, so is $\funcspec $.  
\end{theorem}

\begin{proof}
The functor $\Phi $ induces a $\cat{D}^{\Sigma }$-functor
$\cat{C}^{\Sigma } \xrightarrow{} (\cat{C}')^{\Sigma }$, which takes the
symmetric sequence $(X_{n})$ to the symmetric sequence $(\Phi X_{n})$.  It
follows that $\Phi $ induces a $\spec $-functor $\funcspec
\mathcolon \modspec  \xrightarrow{} \genspec{\cat{C}'}{K}$,
that takes the symmetric spectrum $(X_{n})$ to the symmetric spectrum
$(\Phi X_{n})$, with structure maps 
\[
\Phi X_{n}\otimes K \cong \Phi (X_{n}\otimes K) \xrightarrow{}\Phi X_{n+1}.  
\]
Let $\Gamma $ denote the right adjoint of $\Phi $.  Then the right adjoint of
$\funcspec $ is $Sp^{\Sigma }(\Gamma )$, which takes the symmetric spectrum
$(Y_{n})$ to the symmetric spectrum $(\Gamma Y_{n})$, with structure maps
adjoint to the composite 
\[
\Phi (\Gamma Y_{n}\otimes K) \cong \Phi \Gamma Y_{n}\otimes K
\xrightarrow{} Y_{n}\otimes K \xrightarrow{}Y_{n+1}.
\]
Since $\Ev _{n} Sp^{\Sigma }(\Gamma )=\Gamma \Ev _{n}$, it follows that
$\funcspec F_{n}=F_{n}\Phi $.

It is clear that $Sp^{\Sigma }(\Gamma )$ preserves level fibrations and level
equivalences, so $\funcspec $ is a Quillen functor with respect to the
projective model structure.  In view of Theorem~\ref{thm-localization},
to see that $\funcspec $ defines a Quillen functor with respect to the
stable model structures, it suffices to show that $\funcspec
(s_{n}^{QC})$ is a stable equivalence for all domains and codomains of
$C$ of the generating cofibrations of $\cat{C}$.  But one can readily
verify that $\funcspec (s_{n}^{QC})=s_{n}^{\Phi QC}$, which is a stable
equivalence as required.  Thus $\funcspec $ is a Quillen functor with
respect to the stable model structures. 

If $\Phi $ is a Quillen equivalence, one can easily check that
$\funcspec $ is a Quillen equivalence with respect to the projective
model structure.  To see that it is still a Quillen equivalence with
respect to the stable model structures, we need only show that
$Sp^{\Sigma }(\Gamma )$ reflects stably fibrant objects, in view of
Proposition~\ref{prop-local-quillen}.  But, if $X$ is level fibrant and
$Sp^{\Sigma }(\Gamma )(X)$ is an $\Omega $-spectrum, this means that the
map $\Gamma X_{n}\xrightarrow{}(\Gamma X_{n+1})^{K}\cong \Gamma
(X_{n+1}^{K})$ is a weak equivalence for all $n$.  Since $\Gamma $ reflects
weak equivalences between fibrant objects, this means that $X$ is an
$\Omega $-spectrum, as required.
\end{proof}

Symmetric spectra are also functorial, in a limited sense, in the
cofibrant object $K$.

\begin{theorem}\label{thm-symm-object}
Suppose $f\mathcolon K\xrightarrow{}K'$ is a weak equivalence of
cofibrant objects of $\cat{D}$, and suppose the domains of the
generating cofibrations of $\cat{D}$ and $\cat{C}$ are cofibrant.  Then
$f$ induces a Quillen equivalence $\genspec{\cat{C}}{f} \mathcolon
\modspec \xrightarrow{}\genspec{\cat{C}}{K'}$ which is natural with
respect to $\cat{D}$-Quillen functors of $\cat{C}$.  
\end{theorem}

\begin{proof}
The map $f$ induces a map of commutative monoids $\Sym (K)
\xrightarrow{}\Sym (K')$.  This induces the usual induction and
restriction adjunction 
\[
\genspec{\cat{C}}{f} \mathcolon \modspec \xrightarrow{}
\genspec{\cat{C}}{K'}.
\]
That is, if $X$ is in $\modspec $, then
$\genspec{\cat{C}}{f}(X)=X\otimes _{\Sym (K)} \Sym (K')$.  Restriction
obviously preserves level fibrations and level equivalences, so is a
Quillen functor with respect to the projective model structure.  One can
easily check that $\genspec{\cat{C}}{f}\circ F_{n}=F_{n}$.  It follows
that $\genspec{\cat{C}}{f}(s_{n}^{QC})$ is the map
\[
F_{n+1}(QC\otimes K) \xrightarrow{} F_{n}QC
\]
in $\genspec{\cat{C}}{K'}$.  The weak equivalence $QC\otimes K
\xrightarrow{}QC\otimes K'$ induces a level equivalence
$F_{n+1}(QC\otimes K) \xrightarrow{}F_{n+1}(QC\otimes K')$.  Since the
map $F_{n+1}(QC\otimes K')\xrightarrow{}F_{n}QC$ is a stable equivalence,
so is the given map.  Thus induction is a Quillen functor with respect
to the stable model structures.

We now prove that induction is a Quillen equivalence between the
projective model structures.  The proof of this is similar to the proof
of Theorem~\ref{thm-equiv}.  That is, restriction certainly reflects
level equivalences between level fibrant objects.  It therefore suffices
to show that the map $X \xrightarrow{}X\otimes _{\Sym (K)}\Sym (K')$ is
a level equivalence for all cofibrant $X$.  The argument of
Theorem~\ref{thm-equiv} will prove this without difficulty.  

To show that induction is a Quillen equivalence between the stable model
structures, we need only check that restriction reflects
stably fibrant objects.  This follows from the fact that the map
$Z^{K'}\xrightarrow{}Z^{K}$ is a weak equivalence for all fibrant $Z$.  
\end{proof}

In particular, it does not matter, up to Quillen equivalence, what model
of the simplicial circle one takes in forming the symmetric spectra
of~\cite{hovey-shipley-smith}.  

\section{Comparison of spectra and symmetric spectra}\label{sec-comparison}

In this section, suppose $\cat{D}$ is a left proper cellular symmetric
monoidal model category, $K$ is a cofibrant object of $\cat{D}$, and
$\cat{C}$ is a left proper cellular $\cat{D}$-model category.  Let $G$
denote the left Quillen endofunctor $GX=X\otimes K$ of $\cat{C}$.  Then
we have two different stabilizations of $\cat{C}$, namely $\BF $ and
$\modspec $.  The object of this section is to compare them.  We show
that $\BF $ and $\modspec $ are related by a chain of Quillen
equivalences whenever the cyclic permutation self-map of $K\otimes
K\otimes K$ is homotopic to the identity.  

This is not the ideal theorem; one might hope for a direct Quillen
equivalence rather than a zigzag of Quillen equivalences, and one might
hope for weaker hypotheses, or even no hypotheses.  However, some
hypotheses are necessary, as pointed out to the author by Jeff Smith.
Indeed, the category $\ho \spec $ is symmetric monoidal, and therefore
$\ho \spec (F_{0}S,F_{0}S)$, the self-maps of the unit, form a
commutative monoid.  If we have a chain of Quillen equivalences between
$\genBF{\cat{D}}{K}$ and $\spec $ that preserves the functors $F_{0}$,
then $\ho \genBF{\cat{D}}{K}(F_{0}S,F_{0}S)$ would also have to be a
commutative monoid.  With sufficiently many hypotheses on $\cat{D}$ and
$K$, for example if $\cat{D}$ is the category of simplicial sets and $K$
is a finite simplicial set, we have seen in in Section~\ref{sec-fin-gen}
that this mapping set is the colimit $\colim \ho \cat{D}(K^{\otimes
n},K^{\otimes n})$.  There are certainly examples where this monoid is
not commutative; for example $K$ could be the mod $p$ Moore space, and
then homology calculations show this colimit is not commutative.  In
fact, this monoid will be commutative if and only if the cyclic
permutation map of $K\otimes K\otimes K$ becomes the identity in $\ho
\cat{D}$ after tensoring with sufficiently many copies of $K$.  Hence we
need some hypothesis on the cyclic permutation map.

The heart of our argument is the following theorem.  The argument of
this theorem can be summarized by saying that commuting stabilization
functors are equivalent, and as such, was suggested to the author by
Mike Hopkins in a different context.  Recall that there are two
different ways to tensor with $K$ on $\BF $; the functor $X\mapsto
X\stensor K$ that is a Quillen equivalence but does not involve the
twist map, and the functor $X\mapsto X\otimes K$ that may not be a
Quillen equivalence but does involve the twist map.  

\begin{theorem}\label{prop-symm-asymm}
Suppose that the functor $X\mapsto X\otimes K$ is a Quillen equivalence
of $\BF $, and also that the domains of the generating cofibrations of
$\cat{D}$ are cofibrant.  Then there is a $\cat{D}$-model category
$\cat{E}$ together with $\cat{D}$-Quillen equivalences $\modspec
\xrightarrow{}\cat{E} \xleftarrow{} \BF $.  Furthermore, we have a
natural isomorphism $[\ho \modspec ](F_{0}A,F_{0}B)\cong [\ho \BF ]
(F_{0}A,F_{0}B)$ for $A,B\in \cat{C}$.  
\end{theorem}

\begin{proof}
We take $\cat{E}=\genBF{\modspec }{K}$, where the functor $G\mathcolon
\modspec \xrightarrow{}\modspec $ used to form $\cat{E}$ is defined by
$GX=X\otimes K$ and is part of the $\cat{D}$-model structure of
$\modspec $ (see the comment following Theorem~\ref{thm-symm-monoid}).
This means that the structure map of $G$ involves the twist map.  By
Theorem~\ref{thm-symm-already-equiv}, $F_{0}\mathcolon \modspec
\xrightarrow{}\cat{E}$ is a $\cat{D}$-Quillen equivalence.  On the other
hand, consider $\genspec{\BF }{K}$, where now we use the action of
$\cat{D}$ on $\BF $ that comes from
Theorem~\ref{thm-spectra-are-simplicial}.  As pointed out in
Remark~\ref{rem-subtle-action}, this means that we are using the functor
$X\mapsto X\otimes K$ to form $\genspec{\BF }{K}$, not the functor
$X\mapsto X\stensor K$.  By hypothesis, this functor is already a
Quillen equivalence, so Theorem~\ref{thm-already-equiv} implies that
$F_{0}\mathcolon \BF \xrightarrow{}\genspec{\BF }{K}$ is a
$\cat{D}$-Quillen equivalence.

It remains to prove that $\cat{E}$ and $\genspec{\BF }{K}$ are
isomorphic as model categories.  This is mostly a matter of unwinding
definitions.  An object of $\cat{E}$ is a set $\{Y_{m,n} \}$ of objects
of $\cat{C}$, where $m,n\geq 0$.  There is an action of $\Sigma _{n}$ on
$Y_{m,n}$, and there are $\Sigma _{n}$-equivariant maps $Y_{m,n}\otimes
K\xrightarrow{\nu }Y_{m+1,n}$ and $Y_{m,n}\otimes K\xrightarrow{\rho
}Y_{m,n+1}$.  In addition, the composite $Y_{m,n} \otimes K^{\otimes
p}\xrightarrow{}Y_{m,n+p}$ is $\Sigma _{n}\times \Sigma
_{p}$-equivariant, and there is a compatibility between $\nu $ and $\rho
$, expressed in the commutativity of the following diagram.
\[
\begin{CD}
Y_{m,n}\otimes K\otimes K @>1\otimes T>> Y_{m,n}\otimes K\otimes K @>\rho
\otimes 1>> Y_{m,n+1}\otimes K \\ 
@V\nu \otimes 1VV @. @VV\nu V \\
Y_{m+1,n}\otimes K @= Y_{m+1,n}\otimes K @>>\rho > Y_{m+1,n+1}
\end{CD}
\]
An object $\{Z_{m,n} \}$ of the category
$\genspec{\BF }{K}$ has the same description if we switch
$m$ and $n$.  There is then an isomorphism of categories between
$\cat{E}$ and $\genspec{\BF }{K}$ which simply switches
$m$ and $n$.  The model structures are also the same.  Indeed, they are
both the localization of the evident bigraded projective model structure
with respect to the maps $F_{m,n+1}(QC\otimes K)\xrightarrow{}F_{m,n}QC$
and $F_{m+1,n}(QC\otimes K)\xrightarrow{}F_{m,n}QC$, where $F_{m,n}$ is
left adjoint to the evaluation functor $\Ev _{m,n}$.  

The natural isomorphism $[\ho \modspec ](F_{0},F_{0}B)\cong [\ho \BF ]
(F_{0}A,F_{0}B)$ follows from the fact that the composites 
\[
\cat{C} \xrightarrow{F_{0}} \modspec \xrightarrow{F_{0}} \cat{E}
\]
and 
\[
\cat{C} \xrightarrow{F_{0}} \BF \xrightarrow{F_{0}} \genspec{\BF
}{K}\cong \cat{E}
\]
are equal.  
\end{proof}

In particular, we have calculated $[\ho \BF ](F_{0}A,F_{0}B)$ in
Corollary~\ref{cor-maps-stable}; when both the hypotheses of that
corollary and the hypotheses of Theorem~\ref{prop-symm-asymm} hold, we
get the expected result 
\[
[\ho \modspec ](F_{0}A,F_{0}B)=\colim \ho \cat{C}(A\smash K^{\smash
n},B\smash K^{\smash n})
\]
for cofibrant $A$ and $B$.

Theorem~\ref{prop-symm-asymm} indicates that we should try to prove that
$-\otimes K$ is a Quillen equivalence of $\BF $.  The only reasonable
way to do this is by comparing this functor to $-\stensor K$, which we
know is a Quillen equivalence.  The basic idea of the proof is to
compare $X\otimes K\otimes K$ to $X\stensor K\stensor K$.  Both of these
spectra have the same spaces, and their structure maps differ precisely
by the cyclic permutation self-map of $K\otimes K\otimes K$.  So if we
knew that this map were the identity, they would be the same spectra.
The hope is then that, if we know that the cyclic permutation map is
homotopic to the identity, these two spectra are equivalent in $\ho \BF
$.  One can in fact construct a map of spectra $X\stensor K\stensor
K\xrightarrow{}R(X\otimes K\otimes K)$, where $R$ is a level fibrant
replacement functor and $X$ is cofibrant, by inductively modifying the
identity map.  Unfortunately, the author does not know how to do this
modification in a natural way, so is unable to prove that the derived
functors $L(X\otimes K\otimes K)$ and $L(X\stensor K\stensor K)$ are
equivalent in this way.  

Instead, we will follow a suggestion of Dan Dugger and construct a new
functor $F$ on cofibrant objects of $\BF $ and natural level
equivalences $FX\xrightarrow{}X\stensor K\stensor K$ and
$FX\xrightarrow{}X\otimes K\otimes K$, for cofibrant $X$.  It will follow
immediately that $L(X\stensor K\stensor K)$ is naturally equivalent to
$L(X\otimes K\otimes K)$, and so that $X\mapsto X\otimes K$ is a Quillen
equivalence on $\BF $.  Unfortunately, to make this construction we will
need to make some unpleasant assumptions that ought to be unnecessary.  

\begin{definition}\label{defn-interval}
Given a symmetric monoidal model category $\cat{D}$ whose unit $S$ is
cofibrant, a \emph{unit interval} in $\cat{D}$ is a cylinder object $I$
for $S$ such that there exists a map $H_{I}\mathcolon I\otimes
I\xrightarrow{}I$ satisfying $H_{I}\circ (1\otimes i_{0})=H_{I}\circ
(i_{0}\otimes 1)=i_{0}\pi $ and $H_{I}\circ (1\otimes i_{1})$ is the
identity.  Here $i_{0},i_{1}\mathcolon S\xrightarrow{}I$ and $\pi
\mathcolon I\xrightarrow{}S$ are the structure maps of $I$.  Given a
cofibrant object $K$ of $\cat{D}$, we say that $K$ is \emph{symmetric}
if there is a unit interval $I$ and a homotopy
\[
H\mathcolon K\otimes K\otimes K\otimes I\xrightarrow{}K\otimes K\otimes K
\]
from the cyclic permutation to the identity.  
\end{definition}

Note that $[0,1]$ is a unit interval in the usual model structure on
compactly generated topological spaces, and $\Delta [1]$ is a unit
interval in the category of topological spaces.  Indeed, the required
map $H_{1}\mathcolon \Delta [1]\times \Delta [1]$ takes both of the
nondegenerate $2$-simplices $011\times 001$ and $001\times 011$ to
$001$.  Similarly, the standard unit interval chain complex of abelian
groups is a unit interval in the projective model structure on chain
complexes.  Also, any symmetric monoidal left Quillen functor preserves unit
intervals.  It follows, for example, that the unstable $A^{1}$-model
category of Morel-Voevodksy has a unit interval.  


Our goal, then, is to prove the following theorem.  

\begin{theorem}\label{thm-wrong-tensor}
Suppose $\cat{D}$ is a symmetric monoidal model category with cofibrant
unit $S$, and $\cat{C}$ is a left proper cellular $\cat{D}$-model
category.  Suppose that $K$ is a cofibrant object of $\cat{D}$, and that
either $K$ is itself symmetric or the domains of the generating
cofibrations of $\cat{C}$ are cofibrant and $K$ is weakly equivalent to
a symmetric object of $\cat{D}$. Then the functor $X\mapsto X\otimes K$
is a Quillen equivalence of $\BF $.
\end{theorem}

This theorem is certainly not the best one can do.  For example, by
considering the analogous functors with more than three tensor factors
of $K$, it should be possible to show that the same theorem holds if
there is a left homotopy between some even permutation of $K^{\otimes
n}$ and the identity.  Also, it seems clear that one should only need
the cyclic permutation, or more generally some even permutation, to be
equal to the identity in $\ho \cat{D}$.  That is, we should not need a
specific left homotopy.  But the author does not know how to remove this
hypothesis.  

In any case, we have the following corollary.

\begin{corollary}\label{cor-comparison}
Suppose $\cat{D}$ is a left proper cellular symmetric monoidal model
category whose unit $S$ is cofibrant, and whose generating cofibrations
can be taken to have cofibrant domains.  Suppose $\cat{C}$ is a left
proper $\cat{D}$-model category.  Suppose $K$ is a cofibrant object of
$\cat{D}$, and either that $K$ is itself symmetric, or else that the
domains of the generating cofibrations of $\cat{C}$ are cofibrant and
$K$ is weakly equivalent to a symmetric object of $\cat{D}$.  Then there
is a $\cat{D}$-model category $\cat{E}$ and $\cat{D}$-Quillen
equivalences
\[
\modspec \xrightarrow{}\cat{E} \xleftarrow{} \BF .
\]
\end{corollary}

We will prove Theorem~\ref{thm-wrong-tensor} in a series of lemmas.  

\begin{lemma}\label{lem-general}
Suppose $\cat{D}$ is a symmetric monoidal model category whose unit $S$
is cofibrant.  Suppose we have a square 
\[
\begin{CD}
A @>f>> X \\
@VrVV @VVsV \\
B @>>g> Y
\end{CD}
\]
in a $\cat{D}$-model category $\cat{C}$, where $A$ and $B$ are
cofibrant, and a left homotopy $H\mathcolon A\otimes I\xrightarrow{}Y$
from $gr$ to $sf$, for some unit interval $I$.  Then there is an object
$B'$ of $\cat{C}$, a weak equivalence $B'\xrightarrow{q}B$, a
commutative square
\[
\begin{CD}
A @>f>> X \\
@Vr'VV @VVsV \\
B' @>>g'> Y
\end{CD}
\]
such that $qr'=r$, and a left homotopy $H'\mathcolon B'\otimes
I\xrightarrow{}Y$ between $gq$ and $g'$.  Furthermore, this construction
is natural in an appropriate sense.  
\end{lemma}

Naturality means that, if we have a map of such homotopy commutative
squares that preserves the homotopies, then we get a map of the
resulting commutative squares that preserves the maps $q$ and $H'$.  The
precise statement is complicated, and we leave it to the reader.  

\begin{proof}
We let $B'$ be the mapping cylinder of $r$.  That is, we take $B'$ to be
the pushout in the diagram below.
\[
\begin{CD}
A @>i_{0}>> A\otimes I \\
@VrVV  @VVhV \\
B @>>j> B'
\end{CD}
\]
The map $r'$ is then the composite $hi_{1}$, and the map $g'$ is the map
that is $g$ on $B$ and $H$ on $A\otimes I$.  It follows that
$g'r'=Hi_{1}=sf$, as required.  The map $q\mathcolon B'\xrightarrow{}B$
is defined to be the identity on $B$ and the composite $A\otimes
I\xrightarrow{\pi }A\xrightarrow{r}B$ on $A\otimes I$.  Since $j$ is a
trivial cofibration (as a pushout of $i_{0}$), it follows that $q$ is a
weak equivalence, and it is clear that $qr'=r$.  We must now construct
the homotopy $H'$.  First note that $B'$ is cofibrant, since $B$ is so
and $j$ is a trivial cofibration, and so $B'\otimes I$ is a cylinder
object for $B'$.  In fact, $B'\otimes I$ is the pushout of $A\otimes
I\otimes I$ and $B\otimes I$ over $A\otimes I$.  Define $H'$ to be the
constant homotopy $B\otimes I\xrightarrow{\pi }B\xrightarrow{g}Y$ on
$B\otimes I$ and the homotopy
\[
A\otimes I\otimes I\xrightarrow{1\otimes H_{I}} A\otimes
I\xrightarrow{H} Y
\]
on $A\otimes I\otimes I$.  The fact that $H_{I}\circ (i_{0}\otimes
1)=i_{0}\pi $ guarantees that $H'$ is well-defined, and the other
conditions on $H_{I}$ guarantee that $H'$ is a left homotopy from $qg$
to $g'$.  We leave the naturality of this construction to the reader. 
\end{proof}

We also need the following lemma about the behavior of unit intervals. 

\begin{lemma}\label{lem-intervals}
Suppose $\cat{D}$ is a symmetric monoidal model category whose unit $S$
is cofibrant.  Let $I$ and $I'$ be unit intervals, and define $J$ by the
pushout diagram below.
\[
\begin{CD}
S @>i_{0}>> I' \\
@Vi_{1}VV @VVj_{1}V \\
I @>>j_{0}> J
\end{CD}
\]
Then $J$ is a unit interval.
\end{lemma}

\begin{proof}
The reader is well-advised to draw a picture in the topological or
simplicial case, from which the proof should be clear.  We think of $J$
as the interval whose left half is $I$ and whose right half is $I'$.  In
particular, $J$ is a cylinder object for $S$, where $i'_{0}=j_{0}i_{0}$
and $i'_{1}=j_{1}i_{1}$.  Then, because the tensor product preserves
pushouts, we can think of $J\otimes J$ as a square consisting of a copy
of $I\otimes I$ in the lower left, a copy of $I\otimes I'$ in the upper
left, a copy of $I'\otimes I$ is the lower right, and a copy of
$I'\otimes I'$ in the upper right.  We define the necessary map
$G\mathcolon J\otimes J\xrightarrow{}J$ by defining $G$ on each
subsquare.  On the lower left, we use the composite $I\otimes
I\xrightarrow{H}I\xrightarrow{j_{0}}J$, where $H$ is the homotopy making
$I$ into a unit interval.  Similarly, on the upper right square, we use
the composite $j_{1}H'$.  On the upper left square we use the constant
homotopy $j_{0}(1\otimes \pi )$, and on the lower right square we use
the constant homotopy $j_{0}(\pi \otimes 1)$. We leave it to the reader
to check that this makes $J$ into a unit interval.
\end{proof}

The importance of these two lemmas for spectra is indicated in the following
consequence.  

\begin{lemma}\label{lem-homotopy-comm}
Suppose $\cat{D}$ is a left proper cellular symmetric monoidal model
category with a unit interval $I$, whose unit $S$ is cofibrant.  Let $K$
be a cofibrant object of $\cat{D}$, and let $\cat{C}$ be a left proper
cellular $\cat{D}$-model category.  Suppose $A,B\in \BF $, where $A$ is
cofibrant, and we have maps $f_{n}\mathcolon A_{n}\xrightarrow{}B_{n}$
for all $n$ and a homotopy $H_{n}\mathcolon A_{n}\otimes K\otimes
I\xrightarrow{}B_{n+1}$ from $f_{n+1}\sigma _{A}$ to $\sigma
_{B}(f_{n}\otimes 1)$, where $\sigma _{(-)}$ is the structure map of the
spectrum $(-)$.  Then there is a spectrum $C$, a level equivalence
$C\xrightarrow{h}A$, and a map of spectra $C\xrightarrow{g}B$ such that
$g_{n}$ is homotopic to $f_{n}h_{n}$.  Furthermore, this construction is
natural in an appropriate sense.  
\end{lemma}

Once again, the naturality involves the homotopies $H_{n}$ as well as
the maps $f_{n}$.  We leave the precise statement to the reader.  

\begin{proof}
We define $C_{n}$, $h_{n}$, $g_{n}$ and a homotopy $H'_{n}:C_{n}\otimes
I_{n}\xrightarrow{}B_{n}$ from $g_{n}$ to $f_{n}h_{n}$, where $I_{n}$ is
a unit interval, inductively on $n$,
using Lemma~\ref{lem-general}.  To get started, we take $C_{0}=A_{0}$,
$h_{0}$ to be the identity, $g_{0}$ to be $f_{0}$, and $H'_{0}$ to be
the constant homotopy (with $I_{0}=I$).  For the inductive step, we apply
Lemma~\ref{lem-general} to the diagram 
\[
\begin{CD}
C_{n} \otimes K @>g_{n}\otimes 1>> B_{n}\otimes K \\
@V\sigma (h_{n}\otimes 1)VV @VV\sigma V \\
A_{n+1} @>>f_{n+1}> B_{n+1}
\end{CD}
\]
and the homotopy obtained as follows.  We have a homotopy
\[
C_{n}\otimes K\otimes I_{n} \xrightarrow{1\otimes T} C_{n}\otimes
I_{n}\otimes K \xrightarrow{H'_{n}\otimes 1} B_{n}\otimes K
\xrightarrow{\sigma } B_{n+1}
\]
from $\sigma (f_{n}\otimes 1)(h_{n}\otimes 1)$ to $\sigma (g_{n}\otimes
1)$. On the other hand, we also have the homotopy $H_{n}(h_{n}\otimes
1)$ from $f_{n+1}\sigma (h_{n}\otimes 1)$ to $\sigma (f_{n}\otimes
1)(h_{n}\otimes 1)$.  We can amalgamate these to get a homotopy
$G_{n}\mathcolon C_{n}\otimes K\otimes I_{n+1}\xrightarrow{}B_{n+1}$
from $f_{n+1}\sigma (h_{n}\otimes 1)$ to $\sigma (g_{n}\otimes 1)$, and
$I_{n+1}$ is still a unit interval by Lemma~\ref{lem-intervals}.  Hence
Lemma~\ref{lem-general} gives us an object $C_{n+1}$, a map $\sigma
\mathcolon C_{n}\otimes K\xrightarrow{}C_{n+1}$, and a map
$g_{n+1}\mathcolon C_{n+1}\xrightarrow{}B_{n+1}$ such that
$g_{n+1}\sigma =\sigma (g_{n}\otimes 1)$.  It also gives us a map
$h_{n+1}\mathcolon C_{n+1}\xrightarrow{}A_{n+1}$ such that
$h_{n+1}\sigma =\sigma (h_{n}\otimes 1)$ and a homotopy
$H'_{n+1}\mathcolon C_{n+1}\otimes I_{n+1}\xrightarrow{}B_{n+1}$ from
$f_{n+1}h_{n+1}$ to $g_{n+1}$.  This completes the induction step and
the proof (we leave naturality to the reader). 
\end{proof}

With this lemma in had we can now give the proof of
Theorem~\ref{thm-wrong-tensor}.  

\begin{proof}[Proof of Theorem~\ref{thm-wrong-tensor}] 
We first reduce to the case where $K$ is itself symmetric.  So suppose
the generating cofibrations of $\cat{C}$ have cofibrant domains, and
suppose $K'$ is symmetric and weakly equivalent to $K$; this means there
are weak equivalences
$K\xrightarrow{}RK\xrightarrow{}RK'\xleftarrow{}K'$, where $R$ denotes a
fibrant replacement functor.  This means the total left derived functors
$X\mapsto X\otimes ^{L}K$ and $X\mapsto X\otimes ^{L}K'$ are naturally
isomorphic on the homotopy category of any $\cat{D}$-model category.  In
particular, it suffices to show that $X\mapsto X\otimes K'$ is a Quillen
equivalence on $\genBF{\cat{C}}{K}$.  On the other hand, by
Theorem~\ref{thm-equiv}, there are $\cat{D}$-Quillen equivalences
\[
\genBF{\cat{C}}{K} \xrightarrow{} \genBF{\cat{C}}{RK}
\xrightarrow{}\genBF{\cat{C}}{RK'} \xleftarrow{} \genBF{\cat{C}}{K'}.  
\]
It therefore suffices to show that $X\mapsto X\otimes K'$ is a Quillen
equivalence of $\genBF{\cat{C}}{K'}$; that is, we can assume that $K$
itself is symmetric.  

Let $H$ denote the given homotopy from the cyclic permutation to the
identity of $K\otimes K\otimes K$.  Let $X$ be a cofibrant spectrum, let
$\widetilde{\sigma }$ denote the structure map of $X\stensor K\stensor
K$, and let $\sigma $ denote the structure map of $X\otimes K\otimes K$.
These two structure maps differ by the cyclic permutation, and therefore
we are in the situation of Lemma~\ref{lem-homotopy-comm}, with
$A=X\stensor K\stensor K$, $B=X\otimes K\otimes K$, $f_{n}$ equal to the
identity map, and $H_{n}=(\sigma _{X}\otimes 1\otimes 1)(X_{n+1}\otimes
H)$.  It follows that we get a functor $F$ defined on cofibrant objects
of $\BF $ and natural level equivalences $FX\xrightarrow{h}X\stensor
K\stensor K$ and $FX\xrightarrow{g}X\otimes K\otimes K$, where the
latter map is a level equivalence since $g_{n}$ is homotopic to $h_{n}$.
Thus the total left derived functors of $(-)\stensor K\stensor K$ and
$(-)\otimes K\otimes K$ are naturally isomorphic.  Since we know already
that $(-)\stensor K\stensor K$ is a Quillen equivalence, so is
$(-)\otimes K\otimes K$, and hence so is $(-)\otimes K$.
\end{proof}

\appendix
\section{Cellular model categories}\label{sec-cell}

In this section we define cellular model categories and show that the
projective model structures on $\BF $ and $\modspec $ are cellular when
$\cat{C}$ is so.  This is necessary to be sure that the Bousfield
localizations used in the paper do in fact exist.  The definitions in
this section are taken from~\cite{hirschhorn}.  Throughout this section,
then, $G$ will be a left Quillen endofunctor of $\cat{C}$; when we refer
to $\modspec $, we will be thinking of $\cat{C}$ as a $\cat{D}$-model
category, where $\cat{D}$ is some symmetric monoidal model category, and
of $K$ as a cofibrant object of $\cat{D}$.  

A cellular model category is a special kind of cofibrantly generated
model category.  Three additional hypotheses are needed.  

\begin{definition}\label{defn-cellular}
A model category $\cat{C}$ is \emph{cellular} if there is a set of
cofibrations $I$ and a set of trivial cofibrations $J$ making $\cat{C}$
into a cofibrantly generated model category and also satisfying the
following conditions. 
\begin{enumerate}
\item The domains and codomains of $I$ are compact relative to $I$.  
\item The domains of $J$ are small relative to the cofibrations.
\item Cofibrations are effective monomorphisms.  
\end{enumerate}
\end{definition}

The first hypothesis above requires considerable explanation,
which we will provide below.  We first point out that the second
hypothesis will hold in the projective model structure on $\BF $ or
$\modspec $ when it holds in $\cat{C}$.  

\begin{lemma}\label{lem-domains-small}
Suppose $\cat{C}$ is a cofibrantly generated model category with
generating cofibrations $I$ and generating trivial cofibrations $J$, and
$G$ is a left Quillen endofunctor of $\cat{C}$.  If the domains of $J$
are small relative to the cofibrations in $\cat{C}$, then the domains of
the generating trivial cofibrations $J_{G}$ of the projective model
structure on $\BF $ \ulp $\modspec $\urp{} are small
relative to the cofibrations in $\BF $ \ulp $\modspec $\urp .
\end{lemma}

\begin{proof}
For $\BF $, this follows immediately from the definition of $J_{G}$,
Proposition~\ref{prop-smallness-in-G-spec}, and
Proposition~\ref{prop-proj-cofs-are-level}.  The proof for $\modspec $
is similar.  
\end{proof}

We now discuss the third hypothesis.  

\begin{definition}\label{defn-eff-mono}
Suppose $\cat{C}$ is a category.  A map $f\mathcolon X\xrightarrow{}Y$
is an \emph{effective monomorphism} if $f$ is the equalizer of the two
obvious maps $Y\rightrightarrows Y\amalg _{X}Y$.  
\end{definition}

\begin{proposition}\label{prop-eff-mono}
Suppose $\cat{C}$ is a cofibrantly generated model category and $G$ is a
left Quillen endofunctor of $\cat{C}$.  If cofibrations are effective
monomorphisms in $\cat{C}$, then level cofibrations, and in particular
projective cofibrations, are effective monomorphisms in $\BF $ and in
$\modspec $.  
\end{proposition}

\begin{proof}
This is immediate, since limits in $\BF $ and $\modspec $ are taken
levelwise.  
\end{proof}

We must now define compactness.  This will involve some preliminary
definitions.  These definitions are extremely technical; the reader is
advised to keep the example of CW-complexes in mind. 

\begin{definition}\label{defn-presented-cell}
Suppose $I$ is a set of maps in a cocomplete category.  A
\emph{relative $I$-cell complex} is a map which can be written as the
transfinite composition of pushouts of coproducts of maps of $I$.  That
is, given a relative $I$-cell complex $f$, there is an ordinal $\lambda
$ and a $\lambda $-sequence $X\mathcolon \lambda \xrightarrow{}\cat{C}$
and a collection $\{(T^{\beta },e^{\beta },h^{\beta })_{\beta <\lambda }
\}$ satisfying the following properties.
\begin{enumerate}
\item $f$ is isomorphic to the transfinite composition of $X$. 
\item Each $T^{\beta }$ is a set. 
\item Each $e^{\beta }$ is a function $e^{\beta }\mathcolon T^{\beta
}\xrightarrow{}I$.  
\item Given $\beta <\lambda $ and $i\in T^{\beta }$, if $e^{\beta
}_{i}\mathcolon C_{i}\xrightarrow{}D_{i}$ is the image of $i$ under
$e^{\beta }$, then $h^{\beta }_{i}$ is a map $h^{\beta }_{i}\mathcolon
C_{i}\xrightarrow{}X_{\beta }$.  
\item Each $X_{\beta +1}$ is the pushout in the diagram 
\[
\begin{CD}
\coprod _{T^{\beta }} C_{i} @>>> \coprod _{T^{\beta }} D_{i} \\
@V\coprod h^{\beta }_{i}VV @VVV \\
X_{\beta } @>>> X_{\beta +1}
\end{CD}
\]
\end{enumerate}
The ordinal $\lambda $ together with the $\lambda $-sequence $X$ and the
collection $\{(T^{\beta },e^{\beta },h^{\beta })_{\beta <\lambda } \}$
is called a \emph{presentation} of $f$.  The set $\coprod _{\beta
}T_{\beta }$ is the \emph{set of cells} of $f$, and given a cell $e$,
its \emph{presentation ordinal} is the ordinal $\beta $ such that $e\in 
T^{\beta }$.  The \emph{presentation ordinal} of $f$ is $\lambda $.  
\end{definition}

We also need to define subcomplexes of relative $I$-cell complexes.  

\begin{definition}\label{defn-subcomplex}
Suppoe $\cat{C}$ is a cocomplete category and $I$ is a set of maps in
$\cat{C}$.  Given a presentation $\lambda $, $X\mathcolon \lambda
\xrightarrow{}\cat{C}$, and $\{ (T^{\beta },e^{\beta
},h^{\beta })_{\beta <\lambda } \}$ of a map $f$ as a relative $I$-cell
complex, a \emph{subcomplex} of $f$ (or really of the presentation of
$f$), is a collection $\{ (\widetilde{T}^{\beta },\widetilde{e}^{\beta
},\widetilde{h}^{\beta })_{\beta <\lambda } \}$ such that the following
properties hold.  
\begin{enumerate}
\item Every $\widetilde{T}^{\beta }$ is a subset of $T^{\beta }$, and
$\widetilde{e}^{\beta }$ is the restriction of $e^{\beta }$ to
$\widetilde{T}^{\beta }$.  
\item There is a $\lambda $-sequence $\widetilde{X}\mathcolon \lambda
\xrightarrow{}\cat{C}$ such that $\widetilde{X}_{0}=X_{0}$ and a map of
$\lambda $-sequences $\widetilde{X}\xrightarrow{}X$ such that, for every
$\beta <\lambda $ and $i\in \widetilde{T}^{\beta }$, the map
$\widetilde{h}^{\beta }_{i}\mathcolon
C_{i}\xrightarrow{}\widetilde{X}_{\beta }$ is a factorization of
$h^{\beta }_{i}\mathcolon C_{i}\xrightarrow{}X_{\beta }$ through the map
$\widetilde{X}_{\beta }\xrightarrow{}X_{\beta }$.  
\item Every $X_{\beta +1}$ is the pushout in the diagram 
\[
\begin{CD}
\coprod _{\widetilde{T}^{\beta }} C_{i} @>>> \coprod
_{\widetilde{T}^{\beta }} D_{i} \\
@V\coprod \widetilde{h}^{\beta }_{i}VV @VVV \\
\widetilde{X}_{\beta } @>>> \widetilde{X}_{\beta +1}
\end{CD}
\]
\end{enumerate}
\end{definition}

Given a subcomplex of $f$, the \emph{size} of that subcomplex is the
cardinality of its set of cells $\coprod _{\beta <\lambda
}\widetilde{T}^{\beta }$.  Usually, $I$ will be a set of cofibrations in
a model category where the cofibrations are essential monomorphisms.
This condition guarantees that a subcomplex is uniquely determined by
its set of cells~\cite[Proposition~12.5.9]{hirschhorn}.  Of course,
every set of cells does not give rise to a subcomplex.

We can now define compactness. 

\begin{definition}\label{defn-compact}
Suppose $\cat{C}$ is a cocomplete category and $I$ is a set of maps in
$\cat{C}$.  
\begin{enumerate}
\item Given a cardinal $\kappa $, an object $X$ is \emph{$\kappa
$-compact relative to $I$} if, for every relative $I$-cell complex
$f\mathcolon Y\xrightarrow{}Z$ and for every presentation of $f$, every
map $X\xrightarrow{}Z$ factors through a subcomplex of size at most
$\kappa $.  
\item An object $X$ is \emph{compact relative to $I$} if $X$ is $\kappa
$-compact relative to $I$ for some cardinal $\kappa $.  
\end{enumerate}
\end{definition}

The following proposition is adapted from an argument of Phil
Hirschhorn's.  

\begin{proposition}\label{prop-compactness}
Suppose $\cat{C}$ is a cellular model category with
generating cofibrations $I$, and $G$ is a left Quillen endofunctor of
$\cat{C}$.  Let $A$ be a domain or codomain of $I$.  Then $F_{n}A$ is
compact relative to $I_{G}$ in $\BF $ or $\modspec $.
\end{proposition}

\begin{proof}
We will prove the proposition only for $\BF $, as the $\modspec $ case
is similar.  Throughout this proof we will use
Proposition~\ref{prop-eff-mono}, which guarantees that subcomplexes in
$\cat{C}_{G}$ are determined by their cells.  Choose an infinite
cardinal $\gamma $ such that the domains and codomains of $I$ are all
$\gamma $-compact relative to $I$.  When dealing with relative $I$-cell
complexes, we can assume that we have a presentation as a transfinite
composition of pushouts of maps of $I$, rather than as a transfinite
composition of pushouts of coproducts of maps of $I$,
using~\cite[Lemma~2.1.13]{hovey-model} or~\cite[Section
12.2]{hirschhorn}.  A similar comment holds for relative $I_{G}$-cell
complexes.  We will proceed by transfinite induction on $\beta $, where
the induction hypothesis is that for every presented relative
$I_{G}$-cell complex $f\mathcolon X\xrightarrow{}Y$ whose presentation
ordinal is $\leq \beta $, and for every map $F_{n}A\xrightarrow{f}Y$
where $n$ is an integer and $A$ is a domain or codomain of $I$, $f$
factors through a subcomplex with at most $\gamma $ $I_{G}$-cells.
Getting the induction started is easy.  For the induction step, suppose
the induction hypothesis holds for all ordinals $\alpha <\beta $, and
suppose we have a presentation
\[
X=X^{0}\xrightarrow{}X^{1}\xrightarrow{} \dots \xrightarrow{} X^{\alpha
} \xrightarrow{} \dots X^{\beta }=Y
\]
of $f\mathcolon X\xrightarrow{}Y$ as a transfinite composition of
pushouts of maps of $I_{G}$.  Then the boundary of each $I_{G}$-cell of
this presentation is represented by a map $F_{m}C\xrightarrow{}X^{\alpha
}$, for some $\alpha <\beta $, some $m\geq 0$, and some domain $C$ of a
map of $I$.  This map factors through a subcomplex with at most $\gamma
$ $I_{G}$-cells, by induction.  It follows that the $I_{G}$-cell itself
is contained in a subcomplex of at most $\gamma $ $I_{G}$-cells, since
we can just attach the interior of the $I_{G}$-cell to the given
subcomplex.

Now suppose we have an arbitrary map $F_{n}A\xrightarrow{}Y$, where $A$
is a domain or codomain of a map of $I$.  Such a map is determined by a
map $A\xrightarrow{}Y_{n}$ in $\cat{C}$.  The map $f_{n}\mathcolon
X_{n}\xrightarrow{}Y_{n}$ is the transfinite composition of the
cofibrations $X^{\alpha }_{n}\xrightarrow{}X^{\alpha +1}_{n}$.  For each
$\alpha $, there is an $m$ and a map $h$ of $I$ such that $X^{\alpha
}_{n}\xrightarrow{}X^{\alpha +1}_{n}$ is the pushout of $G^{m}h$, where
we interpret $G^{m}h$ as the identity map if $m$ is negative.
Apply~\cite[Lemma~12.4.19]{hirschhorn} to write the $\beta $-sequence
$X^{\alpha }_{n}$ as a retract of a $\beta $-sequence
\[
X_{n}=Z^{0}\xrightarrow{}Z^{1} \xrightarrow{} \dots \xrightarrow{}
Z^{\alpha } \xrightarrow{} \dots Z^{\beta }=Z
\]
where each map $Z^{\alpha }\xrightarrow{}Z^{\alpha +1}$ is a relative
$I$-cell complex.  We denote the retraction by $r\mathcolon
Z\xrightarrow{}Y_{n}$, noting that the restriction of $r$ to $Z^{\alpha
}$ factors (uniquely) through $X^{\alpha }_{n}$.  We can think of the
entire map $X_{n}\xrightarrow{}Z$ as a relative $I$-cell complex, each
cell $e$ of which appears in the relative $I$-cell complex
$Z^{t(e)}\xrightarrow{}Z^{t(e)+1}$ for some unique ordinal $t(e)$, and
so has associated to it the $I_{G}$-cell $c(e)$ of $f$ used to form
$X^{t(e) }\xrightarrow{}X^{t(e)+1}$.  The composite
$A\xrightarrow{}Y_{n}\xrightarrow{}Z$ then factors through a subcomplex
$V$ with at most $\gamma $ $I$-cells.  The proof will be completed if we can
find a subcomplex $W$ of $Y$ with at most $\gamma $ $I_{G}$-cells such
that the restriction of $r$ to $V$ factors through $W_{n}$.

Take $W$ to be a subcomplex of $Y$ containing the
$I_{G}$-cells $c(e)$ as $e$ runs through the cells of $V$.  Then $W_{n}$
contains $r(e)$ for every cell of $V$, so $W_{n}$ contains $rV$, as
required.  Furthermore, since each $I_{G}$-cell $c(e)$ lies in a
subcomplex with no more that $\gamma $ $I_{G}$-cells, and $V$ has no
more than $\gamma $ cells, there is a choice for $W$ which has no more
than $\gamma ^{2}=\gamma $ cells.  This completes the induction step and
the proof.  
\end{proof}

Altogether then, we have the following theorem.

\begin{theorem}\label{thm-left-proper-cellular}
Suppose $\cat{C}$ is a left proper cellular model category, and $G$ is a
left Quillen endofunctor on $\cat{C}$.  Then the category $\BF $ of
$G$-spectra and the category $\modspec $ of symmetric spectra, with the
projective model structures, are left proper cellular model categories.  
\end{theorem}


\providecommand{\bysame}{\leavevmode\hbox to3em{\hrulefill}\thinspace}

\end{document}